\documentclass[smallcondensed]{svjour3}

\usepackage[utf8]{inputenc}
\usepackage[T1]{fontenc}
\usepackage{amsmath,amssymb,fancybox,pgf,array}

\newcommand{\immagine}[3][8]{
\begin{figure}[htb]
\begin{center}
\includegraphics[width=#1cm]{#2.eps}
\caption{#3}
\label{fig:#2}
\end{center}
\end{figure}}

\renewcommand{\b}[1]{{\bf #1}}

\newcommand{\mmfunz}[5]{
$$ #1 : \left\{ \begin{array}{ccl}
 #2 & \rightarrow & #3 \\
 #4 & \mapsto& #5   \end{array} \right. $$}
\newcommand{\fz}[3]{#1:\, #2 \rightarrow #3}

\renewcommand{\r}[1]{(\ref{#1})}
\newcommand{\bi}{\begin{itemize}}
\newcommand{\ei}{\end{itemize}}
\newcommand{\be}{\begin{enumerate}}
\newcommand{\ee}{\end{enumerate}}

\newcommand{\bd}{\begin{description}}
\newcommand{\ed}{\end{description}}
\renewcommand{\i}{\item}

\newcommand{\bqn}{\begin{eqnarray}}
\newcommand{\eqn}{\end{eqnarray}}
\newcommand{\eqnn}{\nonumber\end{eqnarray}}
\newcommand{\eqnl}[1]{\label{#1}\end{eqnarray}}

\newcommand{\nn}{\nonumber\\}

\newcommand{\ba}[1]{\begin{array}{#1}}
\newcommand{\ea}{\end{array}}

\newcommand{\R}{\mathbb{R}}

\newcommand{\N}{\mathbb{N}}

\newcommand{\nei}{neighborhood}

\newcommand{\bproof}{\begin{proof}}
\newcommand{\eproof}{\end{proof}}

\newcommand{\bt}{\begin{theorem}}
\newcommand{\et}{\end{theorem}}
\newcommand{\bl}{\begin{lemma}}
\newcommand{\el}{\end{lemma}}
\newcommand{\bp}{\begin{proposition}}
\newcommand{\ep}{\end{proposition}}
\newcommand{\bc}{\begin{corollary}}
\newcommand{\ec}{\end{corollary}}
\newcommand{\bdeff}{\begin{definition}}
\newcommand{\edeff}{\end{definition}}
\newcommand{\brem}{\begin{remark}\rm}
\newcommand{\erem}{\end{remark}}
\newcommand{\auth}[1]{{\sc #1}}
\newcommand{\tit}[1]{{\rm #1}}

\newcommand{\jou}[1]{{\it #1}}

\newcommand{\pp}[1]{pp.~#1}
\newcommand{\lam}{\lambda}
\newcommand{\al}{\alpha}
\newcommand{\eps}{\varepsilon}
\newcommand{\ga}{\gamma}

\newcommand{\Id}{\mathrm{Id}}

\newcommand{\Pt}[1]{\left( #1 \right)}
\newcommand{\Pg}[1]{\left\{ #1 \right\}}
\newcommand{\Pq}[1]{\left[ #1 \right] }
\newcommand{\Pa}[1]{\langle #1 \rangle}
\newcommand{\Pabs}[1]{\left| #1 \right|}

\newcommand{\flow}{\Phi}

\newcommand{\dt}{{\Delta t}}

\newcommand{\Prr}{\mathcal{P}^{ac}_c}
\renewcommand{\Pr}{\mathcal{P}_p}
\newenvironment{pschema}[1]{\vspace{3mm}
\noindent\begin{Sbox}\begin{minipage}{\textwidth}\vspace{2mm}\begin{center}{\large \schema{#1}}\vspace{5mm}\\
\begin{minipage}{0.9\textwidth}}{\end{minipage}\end{center}\vspace{2mm}\end{minipage}\end{Sbox}\fbox{\TheSbox}\vspace{3mm}}
\newcommand{\schema}[1]{\b{\sc #1}}

\newcommand{\Wp}[2][p]{W_{#1}\Pt{#2}}
\newcommand{\Pceil}[1]{{\lceil #1 \rceil}}
\newcommand{\Chi}{{\chi}}
\newcommand{\grid}{{
\begin{pgfpicture}{0cm}{0cm}{.3cm}{.3cm}
\pgfline{\pgfxy(0,0)}{\pgfxy(0,.3)}
\pgfline{\pgfxy(.1,0)}{\pgfxy(.1,.3)}
\pgfline{\pgfxy(.2,0)}{\pgfxy(.2,.3)}
\pgfline{\pgfxy(.3,0)}{\pgfxy(.3,.3)}

\pgfline{\pgfxy(0,0)}{\pgfxy(.3,0)}
\pgfline{\pgfxy(0,.1)}{\pgfxy(.3,.1)}
\pgfline{\pgfxy(0,.2)}{\pgfxy(.3,.2)}
\pgfline{\pgfxy(0,.3)}{\pgfxy(.3,.3)}

\end{pgfpicture}
}}
\newcommand{\rettangolo}[4]{
\pgfline{\pgfxy(#1,#2)}{\pgfxy(#1+#3,#2)}
\pgfline{\pgfxy(#1,#2)}{\pgfxy(#1,#2+#4)}
\pgfline{\pgfxy(#1,#2+#4)}{\pgfxy(#1+#3,#2+#4)}
\pgfline{\pgfxy(#1+#3,#2)}{\pgfxy(#1+#3,#2+#4)}
}

\newcommand{\quadrato}[3]{
\rettangolo{#1}{#2}{#3}{#3}}

\newcommand{\h}{{\Delta x}}
\newcommand{\dx}{{\Delta x}}
\begin{document}

\title{Transport equation with nonlocal velocity\\in Wasserstein spaces:\\convergence of numerical schemes}

\titlerunning{Transport equation with nonlocal velocity in Wasserstein spaces}

\author{Benedetto Piccoli\thanks{Department of Mathematical Sciences, Rutgers University - Camden, Camden, NJ. {\tt piccoli@camden.rutgers.edu}}, Francesco Rossi\thanks{Aix-Marseille Univ, LSIS, 13013, Marseille, France. {\tt francesco.rossi@lsis.org}}}

\maketitle

\begin{abstract}
\noindent Motivated by pedestrian modelling, we study evolution of measures in the Wasserstein space. In particular, we consider the Cauchy problem for a transport equation, where the velocity field depends on the measure itself. 

\noindent We deal with numerical schemes for this problem and prove convergence of a Lagrangian scheme to the solution, when the discretization parameters approach zero. We also prove convergence of an Eulerian scheme, under more strict hypotheses. Both schemes are discretizations of the push-forward formula defined by the transport equation. As a by-product, we obtain  existence and uniqueness of the solution.

\noindent All the results of convergence are proved with respect to the Wasserstein distance. We also show that $L^1$ spaces are not natural for such equations, since we lose uniqueness of the solution.

\end{abstract}

\vspace{.3cm}


\vspace{1cm}

Recently, when addressing modelling of pedestrian crowd motions, various authors used measures to represent relevant quantities, such as pedestrian density (see e.g. \cite{ben-multiscale,evers,maury1,maury2,ben-ARMA,ben-domain,tosin}). The density follows an evolution prescribed by a flow map of a velocity function, that is typically composed of two terms: a first one called desired velocity, depending on the geometry of the state space only, and a second one called interaction velocity, depending
on the position of the other pedestrians, thus on the whole measure. The resulting dynamics is a transport equation of the type:
\bqn
\begin{cases}
\partial_t\mu+\nabla\cdot(v \mu)=0\\
\mu_{|_{t=0}}=\mu_0,
\end{cases}
\label{e-cauchy}
\eqn
where $v=v[\mu]$. Since each pedestrian interacts with others in a surrounding area, the term $v$ is nonlocal, often expressed as an operator depending on a compactly supported kernel. Moreover, it is very useful to work in a Wasserstein space, i.e. the space of probability measures endowed with the Wasserstein distance. This distance is defined in terms of the solution to the optimal transportation problem \`{a} la Monge-Kantorovich. Also, numerical schemes deriving from the push-forward of the measure by the discretized flow, proved to be particularly convenient.\\

The aim of the present paper is to study convergence of numerical schemes for \r{e-cauchy} with $v=v[\mu]$. The schemes we study are discretizations of the push-forward formula defined by  \r{e-cauchy}. Our main result in this context is that schemes in Lagrangian form are the more suitable, converging for any choice of space-time discretization. On the other side, Eulerian ones require quite restrictive assumptions for convergence. As a by-product, convergence of the schemes is used to prove existence and uniqueness of the solution for \r{e-cauchy}, as proved in \cite{ambrosio}.

All along the article, we use the following notation: we deal with measures $\mu$ in $\Prr$, the space of probability measures on $\R^n$ {with compact support and absolutely continuous with respect to the Lebesgue measure}. We also deal with the space of Radon probability measures with finite $p$-moment $\Pr$, on which the Wasserstein distance is defined and finite. Observe that $\Prr\subset \Pr$.

The first problem we address is the existence and uniqueness of a solution for \r{e-cauchy} in the case of $v=v\Pq{\mu}$, i.e. $v$ is a vector field depending on the density itself. In this context, we always assume the following hypotheses:

\newcommand{\Hp}{{\bf (H)}}

\vspace{3mm}
\noindent\begin{Sbox}\begin{minipage}{\textwidth}\vspace{2mm}\begin{center}\Hp\vspace{5mm}\\
\begin{minipage}{0.9\textwidth}
The function
\mmfunz{v\Pq{\mu}}{\Pr}{C^{1}(\R^n)\cap L^\infty(\R^n)}{\mu}{v\Pq{\mu}}
satisfies 
\bi
\i $v\Pq{\mu}$ is uniformly Lipschitz and uniformly bounded, i.e. there exist $L$, $M$ not depending on $\mu$, such that for all $\mu\in\Pr, x,y\in\R^n,$
\bqn
\hspace{-5mm}|v\Pq{\mu}(x)-v\Pq{\mu}(y)|\leq L |x-y|\qquad |v\Pq{\mu}(x)|\leq M.
\eqnn

\i $v$ is a Lipschitz function, i.e. there exists $K$ such that 
\bqn\|v\Pq{\mu}-v\Pq{\nu}\|_{\mathrm{C^0}} \leq K \Wp{\mu,\nu}.
\eqnn
\ei
\end{minipage}\end{center}\vspace{2mm}\end{minipage}\end{Sbox}\fbox{\TheSbox}\vspace{3mm}

Under assumption \Hp, we prove the convergence of a Lagrangian semi-discrete in time scheme for \r{e-cauchy}. We then introduce a complete discrete Lagrangian scheme, for which we also prove convergence to the solution of \r{e-cauchy}. Given a space discretization parameter $\dx$ and a time discretization parameter $\dt$, we prove that the error of the approximation $\mu^L$ of the solution $\mu$ on the time interval $\Pq{0,T}$ satisfies $$W_p(\mu^L(t),\mu(t))\leq a\dx+b\dt$$ with $a,b$ given explicitly in Propositions \ref{p-schema3} and \ref{p-schema4}. It is clear that $a,b$ depend on the constants $L,K,M$ in \Hp, the final time $T$, the metric index $p$ and the dimension of the space $n$. As a consequence, we have weak convergence of $\mu^L$ to $\mu$ for $\dt,\dx\to 0$ only.

We then introduce an Eulerian scheme, that has been first proposed in \cite{ben-ARMA} for modelling of pedestrian dynamics. Similar numerical schemes have been used for a long time in computational fluid dynamics, see e.g. \cite{rew1,rew2,rew3}. Nevertheless, to our knowledge, no formal convergence results in Wasserstein distance were studied. In this article, we prove the convergence of the scheme to the solution of \r{e-cauchy} under more restrictive hypotheses than the Lagrangian scheme. Indeed, we prove convergence in $W_1$ for $\dt, \frac{\dx}{\dt}\rightarrow 0$. We also prove convergence in $W_p$ with $p>1$ for $\dt,\dx 2^{(1-1/p)\frac{T}{\dt}}\to 0$. Also in this case, a precise estimate of the error is given in Proposition \ref{p-schema5}.

\brem The hypothesis \Hp\ is quite restrictive. On one side, the Lipschitz dependence on the Wasserstein distance does not include the case of vector fields depending on point-wise values of the measure, provided a good definition of this quantity is possible. See Section \ref{s-L1} for a more detailed discussion.

On the other side, as shown in \cite{ambrosio}, one can prove existence and uniqueness in the same Wasserstein setting, but with weaker hypotheses than \Hp.\erem

It is interesting to observe that all these estimates are given in terms of the Wasserstein distance. This is natural, since in \Hp\ we  impose Lipschitzianity of $v$ with respect to this distance. For this reason, one could be interested in studying the same equation in $L^1$. This idea is studied in Section \ref{s-L1}, where we replace $W_p$ with $L^1$ in \Hp. The surprising result is that we lose uniqueness of the solution of \r{e-cauchy} under these new hypotheses. For this reason, we don't study numerical schemes in this framework.\\

%The difference between Lagrangian and Eulerian schemes can be explained as follows: for a Lagrangian scheme, we fix a ``particle'' $x$ and compute the evolution in time of its position, i.e. its trajectory $x(t)\in\R^n$. In an Eulerian scheme, we instead fix a ``spatial point'' $P\in\R^n$ and compute the evolution in time of the value of measure at this fixed point, i.e. the value $\mu_t(P)$.\\

The application to pedestrian dynamics also explains the choice of the basic assumptions \Hp, namely that we deal with measures with bounded support. The fact that we deal with absolutely continuous measures forbids us to directly use the results of this paper for multi-scale models, in which the population is modelled by a continuous and discrete part in interaction. Nevertheless, under some assumptions (used e.g. in \cite{ben-multiscale}) we can extend the results of this paper to that context. This important issue is briefly presented in Section \ref{s-multiscala} \\

The structure of the paper is the following. In Section \ref{s1-general} we study properties of the transport equation \r{e-cauchy} with $v$ depending on $\mu$. We also recall the definition of Wasserstein distance and its basic properties. In Section \ref{s-scheme1} we first study properties of the Wasserstein distance under the action of flows. We then study the first numerical scheme, that is the semi-discrete in time Lagrangian \schema{scheme 1}. We prove its convergence in $C(\Pq{0,T},\Pr)$. As a by-product, we also prove in Section \ref{s1-esistenza} existence and uniqueness of the solution of \r{e-cauchy}. 

The main section of the article is Section \ref{s-schemi}, in which we introduce the other numerical schemes and prove their convergence to the solution under suitable conditions. We also discuss in Section \ref{s1-generalizations} how to generalize the main results of the paper to a finite family of interacting measures (i.e. several populations), eventually in a multi-scale modelling.

In Section \ref{s-nuclei} we give examples of velocities $v=v\Pq{\mu}$ that have been introduced in the context of modelling of pedestrian flows. For the examples proposed, we check if \Hp\ are satisfied or not.

Finally, in Section \ref{s-L1} we discuss the choice of Wasserstein distance or $L^1$ distance, both in terms of modelling of pedestrian and in terms of mathematical properties. In particular, we show non-uniqueness of the solution of \r{e-cauchy} when replacing $W_p$ with $L^1$ distance in \Hp.

\section{Transport equation and Wasserstein distance}
\label{s1-general}

In this section we study the equation \r{e-cauchy}, in which the velocity $v$ is either given by an autonomous vector field $v(x)$ or a time-dependent vector field $v_t(x)$. The first goal of the section is to properly define the solution of \r{e-cauchy} in the context of measures.

Given a measurable map $\fz{\ga}{\R^n}{\R^n}$, one can define the push-forward $\ga\#\mu$ of a measure $\mu$ as follows: $$(\ga\#\mu)(E)=\mu(\ga^{-1}(E)).$$
The algebra of $(\ga\#\mu)$-measurable sets is thus  $\Pg{E\subset\R^n\ |\ \ga^{-1}(E) \mbox{ is $\mu$-measurable}}$. The definition of push-forward also implies that, if $\nu=\ga\#\mu$, then for all non-negative functions $\phi$ it holds
\bqn
\int_{\R^n} \Pt{\phi\circ \ga}\,d\mu = \int_{\R^n} \phi\,d\nu.
\eqnl{e-push-strano}
See a proof in \cite[p. 4]{villani}.

Given two positive measures $\mu,\nu$ with the same total mass, one can ask if there exists a measurable map $\gamma$ such that $\gamma\#\mu=\nu$. Moreover, among all the admissible maps, one can define a cost for $\gamma$ of the kind $\int c(x,\gamma(x)) \,d\mu$ with $\fz{c}{\R^n\times \R^n}{[0,\infty)}$ and find the $\gamma$ with minimal cost. This is the idea of optimal transportation, a problem first proposed by Monge in 1781. A complete introduction is given in \cite{villani}. A particular case is given by the cost $c(x,y)=|x-y|^p$ with $p\geq 1$. This problem induces a definition of distance between two measures, called the Wasserstein distance:
\bqn
\label{e-wass}
W_p(\mu,\nu)=\inf\Pg{\Pt{\int_{\R^n} |\gamma(x)-x|^p\,d\mu}^{1/p}\ \mid\ \gamma\#\mu=\nu}.
\eqnn
This is indeed a distance, see e.g. \cite[Ch. 7]{villani}. It is important to remark that the infimum is always attained, since we deal with measures that are absolutely continuous with respect to the Lebesgue measure. Another interesting property is that the distance can be estimated by dividing the measures in different parts, i.e.
\bqn
W_p^p(\mu_1+\mu_2,\nu_1+\nu_2)\leq W_p^p(\mu_1,\nu_1)+W_p^p(\mu_2,\nu_2).
\eqnl{e-wassdivisa}
This formula makes sense for positive measures $\mu_i,\nu_i$ such that $\mu_i(\R^n)=\nu_i(\R^n)$. We give two corollaries of this result, that we use in the following. Given $\mu,\nu$ sharing a common mass $\eta$, we have
\bqn
W_p(\mu,\nu)= W_p(\mu-\eta,\nu-\eta).
\eqnl{e-Wpshared}
Given $\mu=\sum_{i=1}^N \mu_i$ and $\nu=\sum_{i=1}^N \nu_i$ such that $\mu_i(\R^n)=\nu_i(\R^n)\neq 0$, one has
\bqn
W_p(\mu,\nu)\leq \sup_{i=1,\ldots,N} W_p(\mu_i,\nu_i) \mu_i(\R^n)^{-1/p}.
\eqnl{e-Wpsup}
This can be proved as follows.  
\bqn
W_p^p(\mu,\nu)&\leq&\sum_{i=1}^N \Pt{W_p^p(\mu_i,\nu_i) \mu_i(\R^n)^{-1}} \mu_i(\R^n) \leq\nn
&\leq& \Pt{\sup_{i=1,\ldots,N} W_p(\mu_i,\nu_i) \mu_i(\R^n)^{-1/p}}^p \, \sum_{i=1}^N \mu_i(\R^n)=\nn
&=& \Pt{\sup_{i=1,\ldots,N} W_p(\mu_i,\nu_i) \mu_i(\R^n)^{-1/p}}^p.
\eqnn

Also recall that $W_p$ distances are ``ordered'', in the sense that 
\bqn
p_1\leq p_2 \mbox{~~~implies~~~}W_{p_1}(\mu,\nu)\leq W_{p_2}(\mu,\nu).
\eqnl{e-ordine} See \cite[7.1.2]{villani}. This has a direct consequence in our context. Take a function $v$ satisfying \Hp\ for a certain $p_1$ with constants $L,M,K$. Then $v$ satisfies \Hp\ for all $p>p_1$ and the same constants $L,M,K$. The converse is not true, since we only have that, given $X$ bounded metric space, it holds
\bqn
p_1\geq p_2 \mbox{~~~implies~~~}W_{p_1}(\mu,\nu)\leq W_{p_2}^{p_2/p_1}(\mu,\nu)\, diam(X)^{1-p_2/p_1}.
\eqnn
Thus, for $p<p_1$ we only have a condition of H\"olderianity. In particular, the strongest condition of the kind \Hp\ is given for $p=1$.

Given a fixed $p\geq 1$, we use the following distance in $C(\Pq{0,T},\Pr)$:
\bqn
d(\mu,\nu)=\sup_{t\in\Pq{0,T}}W_p(\mu_t,\nu_t).
\eqnl{e-d}

One can also give a Lipschitz family of maps $\ga_t$ to define the push-forward $\mu_t=\ga_t\#\mu$. Under suitable assumptions, one can find a connection between the push-forward of a measure and the transport equation.

\bt[{\cite[Thm 5.34]{villani}}] \label{t-misuretrasporto} Let $\Pt{\ga_t}_{t\in\Pq{0,T}}$ be a locally Lipschitz in time family of diffeomorphisms of $\R^n$, with $\ga_0=\Id$. Let $v=v(t,x)$ be the velocity field associated with the trajectories of $\ga_t$. Given $\mu_0\in\Prr$, and setting $\mu_t:=\gamma_t\#\mu_0$, then $\mu=\mu_{\Pq{0,T}}$ is the unique solution of the linear transport equation
\bqn
\begin{cases}
\partial_t\mu_t+\nabla\cdot(v \mu_t)=0\\
\mu_{|_{t=0}}=\mu_0
\end{cases}
\label{e-misuretrasporto}
\eqn
in $C\Pt{[0,T],\Prr}$, where $\Prr$ is endowed with the weak topology.
\et

The solution of the previous equation is to be intended in the weak sense, i.e. for all functions $f\in C^\infty_c\Pt{\Pt{0,T}\times \R^n}$, it holds $\int_{\Pq{0,T}\times\R^n} \Pt{\partial_t f+\nabla f \cdot v} \,d\mu=0$.

% We remark that, since we require the continuity of $\mu$ w.r.t. $t$, then it is equivalent to state that for a.e. $t\in [0,T]$ it holds that, for all $f\in C^\infty_c\Pt{\R^n}$
%\bqn
%\partial_t\Pt{\int_{\R^n} f \,d\mu_t}=\int_{\R^n} \nabla f \cdot v\, d\mu_t.
%\eqnn
%See a proof in \cite[p. 168]{villani}.

A typical example of an application of Theorem \ref{t-misuretrasporto} is the case in which $v$ is a given Lipschitz vector field, and $\ga_t=\flow^v_t$ is the flow of $v$. We recall that $\flow^v_t(x_0)$ is the unique solution at time $t$ of
\bqn
\begin{cases}
\dot x = v(x)\\
x(0)=x_0.
\end{cases}
\eqnn
Then the velocity field associated to $\gamma_t$ is exactly $v$. One can easily pass to time-dependent vector fields $v_t$, assuming that they generate a flow. This is verified when $v$ is measurable with respect to time, uniformly Lipschitz in $x$ and uniformly bounded. All the results stated above still hold in this case. In this context, both approaches (push-forward of measures and transport equations) are equivalent.

\brem Our problem does not fit exactly into the hypotheses of Theorem \ref{t-misuretrasporto}, since we endow $\Prr$ with the $W_p$ metric, and the corresponding topology, instead of weak topology.

Nevertheless, observe that we always deal with a compactly supported measure $\mu_0$, and that the velocity is bounded. Hence, all $\mu_t$ have compact support, contained in $\mathrm{supp}(\mu_0)+B_{Mt}(0)$. Then change the metric of $\R_n$ outside the support of the $\mu_t$ to have $\R^n$ bounded. In this context, $W_p$ metrizes weak convergence (see \cite[Remark 7.13]{villani}), thus weak topology coincide with $W_p$ topology on $\Prr$.
\erem

\section{Semi-discrete in time Lagrangian scheme}
\label{s-scheme1}

In this section we introduce a semi-discrete in time  Lagrangian scheme to solve \r{e-cauchy} and we prove that it generates a Cauchy sequence in $C(\Pq{0,T},\Pr)$. Since $\Pr$ is complete, then we have a candidate solution for \r{e-cauchy}. We then prove that such candidate is in $C(\Pq{0,T},\Prr)$, and  that it is indeed a solution.

We first need some basic estimates for $W_p$ distances under the action of flows. Beside being interesting in themselves, these estimates will be useful in the following, both for the study of existence of solutions for $v=v\Pq{\mu}$ and the convergence of numerical schemes.

\subsection{Wasserstein distance under the action of flows}
\label{s-Wpflow}

In this section, we prove estimates of the distance $W_p$ under the action of flows $\flow^v_t(\cdot)$. We recall that, given a time $t>0$, the flow $\flow^v_t(\cdot)$ is a diffeomorphism of $\R^n$, thus we can see it as a change of coordinates.
\bp
Let $v$ be an autonomous vector field, Lipschitz with constant $L$ and bounded. Let $\mu,\nu\in\Prr$ be two probability measures. Then
\bqn
W_p(\flow^v_t\#\mu,\flow^v_t\#\nu)\leq e^{\frac{p+1}{p}L t} W_p(\mu,\nu)
\eqnl{e-stime2}
and 
\bqn
W_p(\mu,\flow^v_t\#\mu)\leq \|v\|_{C^0} t.
\eqnl{e-stima1}
\ep
\bproof
Consider a map $\gamma$ realizing $W^p_p(\mu,\nu)=\int_{\R^n} |\gamma(x)-x|^p\, d\mu(x)$ with $\gamma\#\mu=\nu$. Define the map $U(y)=\flow^v_t(\gamma(\flow^{v}_{-t}(y)))$ and remark that $U\#(\flow^{v}_{t}\#\mu)=\flow^{v}_{t}\#\nu$. Thus
\bqn
W^p_p(\flow^{v}_{t}\#\mu,\flow^{v}_{t}\#\nu)&\leq&
\int_{\R^n} |U(y)-y|^p\, d(\flow^{v}_{t}\#\mu)(y)=\nn
&=&\int_{\R^n}  |U(\flow^{v}_{t}(x))-\flow^{v}_{t}(x)|^p \,
\mid\frac{\partial y}{\partial x}\mid\, d(\flow^{v}_{t}\#\mu)(\flow^{v}_{t}(x)),
\eqnn
where the last equality is the change of coordinates $y\mapsto x=\flow^{v}_{-t}(y)$. We estimate the Jacobian $\mid\frac{\partial y}{\partial x}\mid\leq e^{L t}$, see \cite{libro-bressan}. Hence,
\bqn
W^p_p(\flow^{v}_{t}\#\mu,\flow^{v}_{t}\#\nu)\leq
e^{L t} \int_{\R^n} |\flow^{v}_{t}(\gamma(x))-\flow^{v}_{t}(x)|^p \, d\mu(x).
\eqnn
Gronwall's lemma gives $|\flow^{v}_{t}(\gamma(x))-\flow^{v}_{t}(x)|\leq e^{L t} |\gamma(x)-x| $, thus
\bqn
W^p_p(\flow^{v}_{t}\#\mu,\flow^{v}_{t}\#\nu)\leq 
e^{L t}e^{p Lt}\int_{\R^n} |\gamma(x)-x|^p\, d\mu (x) = e^{(p+1)L t} W_p^p(\mu,\nu).
\eqnn

The second estimate is similar. Take $\gamma(x):=\flow^v_t(x)$ and observe that $|\gamma(x)-x|\leq \|v\|_{C^0} t$. Since $\gamma\#\mu=\flow^v_t\#\mu$ by definition, then
\bqn
W_p^p(\mu,\flow^v_t\#\mu)\leq \int_{\R^n} |\gamma(x)-x|^p\, d\mu(x)\leq \|v\|_{C^0}^p t^p \int_{\R^n} d\mu(x)=\|v\|_{C^0}^p t^p.
\eqnn
\eproof

%\brem Also in this case, the results still hold for time-dependent velocities $v_t$, when they are uniformly Lipschitz in $x$, uniformly bounded, continuous in $t$.
%\erem

We now turn our attention to an estimate in which the flows are given by two distinct vector fields.
\bp
\label{t-stime4}
Let $v,w$ be two vector fields, both Lipschitz with constant $L$ and bounded. Let $\mu,\nu\in\Prr$ be two probability measures. Then
\bqn
W_p(\flow^{v}_{t}\#\mu,\flow^{w}_{t}\#\nu)\leq e^{\frac{p+1}{p} L t} W_p(\mu,\nu)+\frac{e^{L t/p}(e^{Lt}-1)}{L} \|v-w\|_{C^0}.
\eqnl{e-stime4}
\ep
\bproof
The proof is similar to the previous one. Consider a map $\gamma$ realizing $W^p_p(\mu,\nu)=\int_{\R^n} |\gamma(x)-x|^p\, d\mu(x)$ with $\gamma\#\mu=\nu$. Define the map $U(y)=\flow^w_t(\gamma(\flow^{v}_{-t}(y)))$ and remark that $U\#(\flow^{v}_{t}\#\mu)=\flow^{w}_{t}\#\nu$. Applying the change of coordinates $y\mapsto x=\flow^{v}_{-t}(y)$, and following the previous proof, we find
\bqn
&&W^p_p(\flow^{v}_{t}\#\mu,\flow^{w}_{t}\#\nu)\leq
e^{L t} \int_{\R^n} |\flow^{w}_{t}(\gamma(x))-\flow^{v}_{t}(x)|^p \, d\mu(x)\leq\nn
&&\leq e^{L t} \int_{\R^n} |\flow^{w}_{t}(\gamma(x))-\flow^{w}_{t}(x)|^p\, d\mu(x)+e^{L t} \int_{\R^n}|\flow^{w}_{t}(x)-\flow^{v}_{t}(x)|^p \, d\mu(x).
\eqnn

Gronwall's lemma gives $|\flow^{w}_{t}(\gamma(x))-\flow^{w}_{t}(x)|\leq e^{L t} |\gamma(x)-x| $. Define $r_\tau(x)=|\flow^{w}_{\tau}(x)-\flow^v_{\tau}(x)|$ and observe that $r_0(x)=0$ and 
\bqn
\dot r_\tau (x)&\leq& |w(\flow^{w}_{\tau}(x))-v(\flow^v_{\tau}(x))|\leq\nn
&\leq& |w(\flow^{w}_{\tau}(x))-w(\flow^v_{\tau}(x))|+|w(\flow^v_{\tau}(x))-v(\flow^v_{\tau}(x))|\leq L r_\tau + \|v-w\|_{C^0}.
\eqnn
thus, applying again Gronwall's lemma, we find $r_t(x)\leq \frac{\|v-w\|_{C^0}}L (e^{Lt}-1)$. Hence
\bqn
W^p_p(\flow^{v}_{t}\#\mu,\flow^{w}_{t}\#\nu)&\leq& e^{(p+1) L t} \int_{\R^n}  |\gamma(x)-x|^p\, d\mu(x)+\nn
&&\hspace{1cm} +e^{L t} \frac{\|v-w\|^p_{C^0}}{L^p} (e^{Lt}-1)^p \int_{\R^n} d\mu(x).
\eqnn
We thus have $W_p(\flow^{w}_{t}\#\mu,\flow^{v}_{t}\#\nu)\leq e^{\frac{p+1}{p} L t} W_p(\mu,\nu)+\frac{e^{L t/p}(e^{Lt}-1)}L \|v-w\|_{C^0}$.
\eproof

\brem These results can be generalized to non-autonomous vector fields $w_t,v_t$, if they generate smooth flows. As already recalled, it is verified if they are measurable with respect to time and uniformly bounded and Lipschitz with respect to space. In this case, in Proposition \ref{t-stime4} we have to replace $\|v-w\|_{C^0}$ with $\sup_t \|v_t-w_t\|_{C^0}$.
\erem

\brem These results can be easily adapted to positive measures that are not probability measures but have finite mass, i.e. $\mu(\R^n)=C\neq 1$. In this case, we have
\bqn
W_p(\flow^{v}_{t}\#\mu,\flow^{w}_{t}\#\nu)\leq e^{\frac{p+1}{p} L t} W_p(\mu,\nu)+\mu(\R^n)^{1/p} \frac{e^{L t/p}(e^{Lt}-1)}L \|v-w\|_{C^0}.
\eqnl{e-stime4massa}
Remark that the formula is symmetric. Indeed, $\mu(\R^n)=\nu(\R^n)$ to have existence of the Wasserstein distance $W_p(\mu,\nu)$.
\erem

\subsection{Definition of the semi-discrete in time Lagrangian scheme}

In this section, we precisely define the the semi-discrete in time Lagrangian scheme and we prove that it defines a Cauchy sequence in $C(\Pq{0,T},\Pr)$. The idea of the scheme is the following : divide the time interval $\Pq{0,T}$ in intervals $\Pq{j\dt,(j+1)\dt}$. For each interval, compute the velocity at the initial time $v_{j\dt}=v\Pq{\mu_{j\dt}}$ and use it as a constant on the whole interval, i.e. compute $\mu_t=\flow^{v_{j\dt}}_{(t-j\dt)}\#\mu_{j\dt}$. 

\begin{pschema}{SCHEME 1\\Lagrangian, semi-discrete in time, exact velocity}
\b{INITIALIZATION:} Fix a time discretization parameter $\dt$.

 Take the starting measure $\mu_0$.\\

%\b{STEP 1:} Given the starting measure $\mu_0$, define $v_0$ on the interval $\Pq{0,\dt}$ as $v_0:=v\Pq{\mu_0}$ and compute the corresponding flow $\flow^{v_0}_t$. 

%For $t\in\Pq{0,\dt}$, define $\mu_t:=\flow^{v_0}_t\#\mu_0$.\\

\b{STEP i+1:} Given $\mu_{i\dt}$, define $v_{i\dt}$ on the interval $\Pq{i\dt,(i+1)\dt}$ as $v_{i\dt}:=v\Pq{\mu_{i\dt}}$ and compute the corresponding flow $\flow^{v_{i\dt}}_t$. 

For $t\in\Pq{i\dt,(i+1)\dt}$, define $\mu_t:=\flow^{v_{i\dt}}_t\#\mu_{i\dt}$.\\

\b{STOP:} when reaching $T$. 
\end{pschema}

Fix a natural number $k$ and divide $\Pq{0,T}$ in $2^k$ intervals, i.e. choose $\dt=\frac{T}{2^k}$. Call $\mu_t^k$ the solution of this numerical scheme. We study the convergence of the sequence $\Pg{\mu_{\Pq{0,T}}^k}_{k\in\N}$ in $C(\Pq{0,T},\Pr)$. We prove the following result.
\bp
\label{p-schema1}
Let $v$ satisfy \Hp, and $\mu_0$ be given. Let $\mu^k=\mu_{\Pq{0,T}}^k$ be constructed by \schema{Scheme 1} with $\dt=\frac{T}{2^k}$. Then, the sequence $\Pg{\mu_{\Pq{0,T}}^k}_{k\in\N}$ is a Cauchy sequence in $C(\Pq{0,T},\Pr)$.
\ep
\bproof
To simplify the notation, we assume that $T=1$. We also estimate independently on $p\geq 1$, and for sufficiently big $k$. We call $m^k_j:=\mu^k_{\frac{j}{2^k}}$ and $v^k_j:=v\Pq{m^k_j}$. The corresponding flow is denoted by $f^{j,k}_t:=\flow^{{v^k_j}}_t$. Remark that we have $\|v^k_j-v^l_i\|_{C^0}\leq K W_p(m^k_j,m^l_i)$.

We  estimate the distance between two successive approximations, that is $d(\mu^k,\mu^{k+1})$. Fix $k\in\N$, $t\in\Pq{0,1}$ and estimate $W_p(\mu^k_t,\mu^{k+1}_t)$. Take $j\in\Pg{0,\ldots,2^{k}-1}$ such that $t\in\Pq{\frac{j}{2^k},\frac{j+1}{2^{k}}}$. We estimate $W_p(\mu^k_t,\mu^{k+1}_t)$ in terms of $W_p(\mu^k_{\frac{j}{2^k}},\mu^{k+1}_{\frac{j}{2^k}})=W_p(m^k_j,m^{k+1}_{2j}),$ that is the distance at the initial time of the interval for the $k$ approximation. We have two cases:
\bi
\i $t\in\Pq{\frac{j}{2^k},\frac{2j+1}{2^{k+1}}}$. Call $t'=t-\frac{j}{2^k}\leq 2^{-k-1}$. We have $\mu^k_t=f^{j,k}_{t'}\# m^k_j$, $\mu^{k+1}_t=f^{2j,k+1}_{t'}\# m^{k+1}_{2j}$. We apply \r{e-stime4} and get
\bqn
W_p(\mu^k_t,\mu^{k+1}_t)&=&W_p(f^{j,k}_{t'}\# m^k_j,f^{2j,k+1}_{t'}\# m^{k+1}_{2j})\leq \nn
&\leq & e^{\frac{p+1}p L t'}W_p(m^k_j,m^{k+1}_{2j})+\frac{e^{L t'/p}(e^{Lt'}-1)}L \| v^k_j-v^{k+1}_{2j}\|_{C^0}\leq\nn
&\leq & \Pt{1+2^{-k}(2L+2K)} W_p(m^k_j,m^{k+1}_{2j}).
\eqnn

\i $t\in\Pq{\frac{2j+1}{2^{k+1}},\frac{j+1}{2^{k}}}$. Call $t'=t-\frac{2j+1}{2^{k+1}}\leq 2^{-k-1}$. We have $\mu^k_t=f^{j,k}_{t'}\# (f^{j,k}_{\frac{1}{2^{k+1}}} \# m^k_j)$, $\mu^{k+1}_t=f^{2j+1,k+1}_{t'}\# (f^{2j,k+1}_{\frac{1}{2^{k+1}}}\# m^{k+1}_{2j})$. The key difference in these two expressions is that, at step $k$ we apply the same $v^k_j$, first for $\frac{1}{2^{k+1}}$ then for $t'$. At step $k+1$, we first apply $v_{2j}^{k+1}$ for $\frac{1}{2^{k+1}}$, then $v_{2j+1}^{k+1}$ for ${t'}$.

We apply \r{e-stime4} and have
\bqn
&&W_p(\mu^k_t,\mu^{k+1}_t)=W_p(f^{j,k}_{t'}\# (f^{j,k}_{\frac{1}{2^{k+1}}}\# m^k_j),f^{2j+1,k+1}_{t'}\# (f^{2j,k+1}_{\frac{1}{2^{k+1}}}\# m^{k+1}_{2j}))\leq\nn
&&\leq \Pt{1+2^{-k+1}L}W_p(f^{j,k}_{\frac{1}{2^{k+1}}}\# m^k_j,f^{2j,k+1}_{\frac{1}{2^{k+1}}}\# m^{k+1}_{2j})+2^{-k+1} \| v^k_j-v^{k+1}_{2j+1}\|_{C^0}.
\eqnn
We estimate the first term using the previous case at time $t=\frac{2j+1}{2^{k+1}}$. For the second term, we have
\bqn
\| v^k_j-v^{k+1}_{2j+1}\|_{C^0}&\leq& K W_p(m^k_j,m^{k+1}_{2j+1})=K W_p(m^k_j,f^{2j,k+1}_{\frac{1}{2^{k+1}}}\# m^{k+1}_{2j})\leq\nn
&\leq& K W_p(m^k_j, m^{k+1}_{2j} )+ K W_p(m^{k+1}_{2j}, f^{2j,k+1}_{\frac{1}{2^{k+1}}}\# m^{k+1}_{2j})\leq\nn
&\leq& K W_p(m^k_j, m^{k+1}_{2j} )+ \frac{K M}{2^{k+1}},
\eqnn
where the last inequality is given by \r{e-stima1}. We finally have 
\bqn
W_p(\mu^k_t,\mu^{k+1}_t)&\leq & (1+2^{-k}C_1) W_p(m^k_j,m^{k+1}_{2j})+2^{-2k} C_2,
\eqnl{e-croce}
with $C_1:=4L+4K+1$, $C_2:=KM$.
\ei
Remark that we both estimates are in terms of $W_p(m^k_j,m^{k+1}_{2j})$, and that we can use \r{e-croce} to estimate on the whole interval $t\in\Pq{\frac{j}{2^k},\frac{j+1}{2^k}}$. Using it recursively in $j=0,\ldots,2^k$, and recalling that $W_p(\mu^k_0,\mu^{k+1}_0)=0$, we have that
\bqn
W_p(\mu^k_t,\mu^{k+1}_t)&\leq & 2^{-2k} C_2\frac{((1+2^{-k}C_1)^{2^k}-1}{(1+2^{-k}C_1)-1}\leq 2^{-k} C_2\frac{e^{C_1}-1}{C_1}.
\eqnn
Since the estimate holds for any $t\in\Pq{0,1}$, we use it to estimate $d(\mu^k,\mu^{k+1})\leq 2^{-k} C_2\frac{e^{C_1}-1}{C_1}$. Since the right-hand side is a Cauchy sequence in $k$, then $\mu^k$ is a Cauchy sequence with respect to $d$.
\eproof

We also prove the continuous dependence of the approximate solution computed via the scheme with respect to the initial data.
\bp
Let $v$ satisfy \Hp, and $\mu_0,\nu_0$ be given. Let $\mu^k$ and $\nu^k$ be constructed using \schema{Scheme 1} with $\dt=\frac{T}{2^k}$ starting from $\mu_0$ and $\nu_0$, respectively.

Then, for a sufficiently small $\dt$, we have
\bqn
d(\mu^k,\nu^k)\leq e^{(4L+4K)T}W_p(\mu_0,\nu_0).
\eqnl{e-stime1}
\ep
\bproof
The proof is similar to the previous one. Fix the following notation $m^k_j:=\mu^k_{j\dt}$, $n^k_j:=\nu^k_{j\dt}$. Then $m^k_{j+1}=\Phi^{v\Pq{m^k_{j}}}_{\dt}\# m^k_j$ and $n^k_{j+1}=\Phi^{v\Pq{n^k_{j}}}_{\dt}\# n^k_j$.

Take $t\in\Pq{j\dt,(j+1)\dt}$ and call $t'=t-j\dt\leq \dt$. Then, using \r{e-stime4} and \Hp, we have 
\bqn
W_p(\mu^k_t,\nu^k_t)&\leq& e^{\frac{p+1}{p}L t'} W_p(m^k_j,n^k_j) + e^{Lt'/p}\Pt{e^{Lt'}-1}\frac{K}{L} W_p(m^k_j,n^k_j)\leq\nn
&\leq& (1+2^{-k}(4L+4K)T)W_p(m^k_j,n^k_j).
\eqnn
Using it recursively in $j$, we have
$$d(\mu^k,\nu^k)\leq (1+2^{-k}(4L+4K)T)^{2^k} W_p(m^k_0,n^k_0)\leq e^{(4L+4K)T}W_p(\mu_0,\nu_0).$$
\eproof

\subsection{Existence and uniqueness of solution with velocity depending on the measure}
\label{s1-esistenza}

In this section we prove existence and uniqueness of a solution for \r{e-cauchy}. The proof is similar to the one given in \cite{ambrosio}. We first show that, if $v$ depends on the measure $\mu_t$ and satisfies \Hp, then there exists a solution of \r{e-cauchy}. The key idea is to use the convergence of the semi-discrete \schema{Scheme 1} and to prove that the limit is indeed a solution of \r{e-cauchy}. We then prove uniqueness of such solution.

We first recall that $\Pr$ is complete with respect to $W_p$, see e.g. \cite[Thm 6.18]{old-new}. As a consequence, $C(\Pq{0,T},\Pr)$ is complete with respect to the distance $d$ defined in \r{e-d}. We now prove that the limit of the semi-discrete \schema{Scheme 1} is a solution of \r{e-cauchy}.

\bt
\label{t-esistenza}
Let $v$ satisfy \Hp, and $\mu_0$ given. Let $\mu^k=\mu^k_{\Pq{0,T}}\in C(\Pq{0,T},\Pr)$ be computed using \schema{Scheme 1} with $\dt=\frac{T}{2^k}$, starting with $\mu^k_0=\mu_0$. Then the limit $\bar{\mu}=\lim_k \mu^k$ exists and is a solution of \r{e-cauchy}.
\et
\bproof
\newcommand{\vt}{\bar{v}_t}
\newcommand{\vkj}{v^k_j}
First remark that $\bar{\mu}$ exists, since $\mu^k$ is a Cauchy sequence, by Proposition \ref{p-schema1}, and $C(\Pq{0,T},\Pr)$ is a complete metric space. It is clear that $\bar{\mu}_0=\mu_0$. To simplify the notation, we assume $T=1$. We also use the notation $\Pa{\mu,f}:=\int_{\R^n} d\mu\, f$ for space integration, and denote $\vt:=v\Pq{\bar{\mu}_t}$, $\vkj:=v\Pq{\mu^k_{j\dt}}$. 

We now prove that, given a test function $f\in C^\infty_c(\Pt{0,1}\times {\R^n})$, it holds
\bqn
\int_0^1 dt\, \Pa{\bar{\mu}_t,\partial_t f+ \vt \cdot \nabla f}=0.
\eqnl{e-musol}
We already know that $\mu^k$ computed via \schema{Scheme 1} satisfies the following
\bqn
\sum_{j=0}^{2^k-1} \int_{j\dt}^{(j+1)\dt}  dt \,\Pa{\mu^k_t,\partial_t f+ \vkj \cdot \nabla f}=0.
\eqnn
Remark that in this case $v$ is evaluated at $t=j\dt$ only. To prove \r{e-musol}, we prove that
\bqn
\lim_\dt \Pabs{\sum_{j=0}^{2^k-1} \int_{j\dt}^{(j+1)\dt} dt \Pq{ \Pa{\bar{\mu}_t ,\partial_t f+ \vt \cdot \nabla f}-\Pa{\mu^k_t,\partial_t f+ \vkj \cdot \nabla f} }}
\eqnl{stima}
is equal to 0. We estimate
\bqn
\r{stima}&\leq& \lim_k \sum_{j=0}^{2^k-1} \int_{j\dt}^{(j+1)\dt} dt \Pabs{\Pa{\bar{\mu}_t,\partial_t f+ \vt \cdot \nabla f}-\Pa{\mu^k_t,\partial_t f+ \vkj \cdot \nabla f}}\leq \nn
&\leq & \lim_k \sum_{j=0}^{2^k-1} \int_{j\dt}^{(j+1)\dt} dt \Pt{\Pabs{\Pa{\bar{\mu}_t-\mu^k_t, \partial_t f}}+\Pabs{\Pa{\bar{\mu}_t, \vt \cdot \nabla f}-\Pa{\mu^k_t, \vkj \cdot \nabla f}}}.
\eqnn

For the first term, define $C_1:=\sup_{[0,T]\times \R^n} |\nabla \partial_t f|$. If $C_1=0$, then $f$ is constant, hence $f=0$ and \r{stima}$=0$. Otherwise, $0<C_1 < \infty$, then $\frac{\partial_t f}{C_1}$ is Lipschitz with respect to space with constant at most $1$, thus $\Pa{ \bar{\mu}_t-\mu^k_t,\frac{\partial_t f}{C_1}}\leq W_1(\bar{\mu}_t,\mu^k_t)$, using the Kantorovich-Rubinstein duality formula, see e.g. \cite[Ch. 1]{villani}. Then recall that Wasserstein distances are ordered, in particular $W_1(\bar{\mu}_t,\mu^k_t)\leq W_p(\bar{\mu}_t,\mu^k_t)$. We now estimate 
\bqn
&&\Pabs{\Pa{\bar{\mu}_t, \vt \cdot \nabla f}-\Pa{\mu^k_t, \vkj \cdot \nabla f}}\leq
\Pabs{\Pa{\bar{\mu}_t, \vt \cdot \nabla f}-\Pa{\bar{\mu}_t, v\Pq{\mu^k_t} \cdot \nabla f}}+\nn
&&+\Pabs{\Pa{\bar{\mu}_t, v\Pq{\mu^k_t} \cdot \nabla f}-\Pa{\mu^k_t, v\Pq{\mu^k_t} \cdot \nabla f}}+
\Pabs{\Pa{\bar{\mu}_t, v\Pq{\mu^k_t} \cdot \nabla f}-\Pa{\mu^k_t, \vkj \cdot \nabla f}}\leq\nn
&\leq & \Pa{\bar{\mu}_t, \Pabs{\vt-v\Pq{\mu^k_t}} \cdot \Pabs{\nabla f}}+\Pabs{\Pa{\bar{\mu}_t-\mu^k_t, v\Pq{\mu^k_t} \cdot \nabla f}}+
\Pa{\bar{\mu}_t, \Pabs{v\Pq{\mu^k_t} - \vkj} \cdot \nabla f}.
\eqnn
Now recall that $\|v\Pq{\mu}-v\Pq{\nu}\|_{C^0}\leq K W_p(\mu,\nu)$. Moreover, both $v$ and $\nabla f$ are Lipschitz and bounded, thus $v\Pq{\mu^k_t} \cdot \nabla f$ is Lipschitz with a certain constant $C_2$ not depending on $t$ or $k$. Thus, again by the Kantorovich-Rubinstein duality formula and the fact that Wasserstein distances are ordered, we have 
\bqn
\Pabs{\Pa{\bar{\mu}_t, \vt \cdot \nabla f}-\Pa{\mu^k_t, \vkj \cdot \nabla f}}\leq (C_3+C_2) W_p(\bar{\mu}_t,\mu^k_t)+ C_3 W_p (\mu^k_t,\mu^k_{j\dt}),
\eqnn
with $C_3:=K\sup_{[0,T]\times \R^n}\Pabs{\nabla f}$. Going back to \r{stima}, we have
\bqn
\r{stima}\leq\lim_\dt \sum_{j=0}^{2^k-1} \int_{j\dt}^{(j+1)\dt} dt \left(W_p(\bar{\mu}_t,\mu^k_t) \Pt{C_1+C_2+ C_3}+ C_3 W_p (\mu^k_t,\mu^k_{{j\dt}}) \right).
\eqnn
We estimate the first term by passing to the supremum over $t$, recalling that $d(\mu,\nu)=\sup_t W_p(\mu_t,\nu_t)$. The last term can be estimated by recalling that $$W_p (\mu^k_t,\mu^k_{{j\dt}})=W_p\Pt{\flow^{v\Pq{\mu^k_{j\dt}}}_{t-j\dt}\# \mu^k_{j\dt},\mu^k_{j\dt}}\leq M(t-j\dt).$$ An integration in time gives 
\bqn
\r{stima}&\leq& \lim_\dt \Pt{d(\bar{\mu},\mu^k) \Pt{C_1+C_2+ C_3} + C_3 M \frac1{2^{k+1}} }=0.
\eqnn

We finally prove that $\bar\mu\in C(\Pq{0,T},\Prr)$, i.e. that $\bar\mu_t\in \Prr$ for all $t$. Define the non-autonomous vector field $v_t:=v\Pq{\bar\mu_t}$ and observe that it is continuous with respect to time, uniformly Lipschitz in space and uniformly bounded. Hence it generates a flow $\Phi_t^{v_t}$, thus a unique solution of \r{e-cauchy}, that is $\Phi_t^{v_t}\# \mu_0\in\Prr$. For uniqueness of the solution, we have $\bar\mu_t=\Phi_t^{v_t}\# \mu_0$. \eproof

We now prove that $\mu_t$ is Lipschitz with respect to time.
\bp
Let $v$ satisfy \Hp, and $\mu$ be a solution of \r{e-cauchy}. Then $\mu$ is Lipschitz with respect to time, i.e.
\bqn
W_p(\mu_t,\mu_s)\leq M |t-s|.
\eqnl{e-lipt}
\ep
\bproof
Since $\mu_t$ is continuous with respect to time, and $v$ is Lipschitz with respect to $\mu$, then $v\Pq{\mu_t}$ is continuous with respect to time. Call $w_t:=v[\mu_t]$ and observe that it is a non-autonomous vector field, continuous with respect to time. Hence the corresponding solution $\nu_t$ of \r{e-cauchy} satisfies \r{e-stima1}, that is equivalent to $W_p(\nu_t,\nu_s)\leq M |t-s|$. By uniqueness of the solution of \r{e-cauchy} for non-autonomous vector fields, the solution is $\nu_t=\mu_t$.
\eproof

\newcommand{\costantona}[1]{{\Pt{2e}^{\Pceil{#1 \max\Pg{2L,8K}}}}}
We now prove the uniqueness of the solution, as well as the continuous dependence on the initial data.
\bt
\label{t-unicita}
Let $\mu,\nu$ be two solutions of \r{e-cauchy} in $C([0,T],\Prr)$, with $v$ satisfying \Hp. Then 
\bqn
W_p(\mu_t,\nu_t)\leq \Pt{2e}^{\Pceil{t \max\Pg{2L,8K}}}W_p(\mu_0,\nu_0),
\eqnl{e-stimasol}
where $\Pceil{x}$ is the ceiling function of $x$, i.e. the smallest integer not less than $x$.

In particular, if $\mu_0=\nu_0$, then $\mu_t=\nu_t$ for all  $t\in[0,T]$.
\et
\bproof
We first observe that $ W_p(\mu_t,\nu_t)$ is Lipschitz with respect to time. Indeed, given two times $t,s\in\Pq{0,T}$, we have
\bqn  W_p(\mu_t,\nu_t)&\leq& W_p(\mu_t,\mu_s)+W_p(\mu_s,\nu_s)+W_p(\nu_s,\nu_t),\eqnn
hence $|W_p(\mu_t,\nu_t)-W_p(\mu_s,\nu_s)|\leq W_p(\mu_t,\mu_s)+W_p(\nu_s,\nu_t)\leq 2M|t-s|$, using \r{e-lipt}.

We now prove the continuous dependence on the initial data. Define two non-autonomous vector fields $f(t,x):=v\Pq{\mu_t}(x)$, $g(t,x):= v\Pq{\nu_t}(x)$. Since $v$ is Lipschitz with respect to $\mu$ and $\mu$ is Lipschitz with respect to time, then both $f$ and $g$ are Lipschitz with respect to time and space, thus they generate a flow. Applying \r{e-stime4} and estimating independently on $p\geq 1$, we have
\bqn
W_p(\mu_\tau,\nu_\tau)\leq e^{ 2L \tau} W_p(\mu_0,\nu_0) + \frac{e^{L\tau}(e^{L\tau}-1)}{L}\sup_{\tau'\in\Pq{0,\tau}}\|v[\mu_{\tau'}]-v[\nu_{\tau'}]\|_{C^0}
\eqnl{e-mio}
Define $\phi(t):=\sup_{\tau\in[0,t]} W_p(\mu_\tau,\nu_\tau)$. Passing to the supremum over $\tau\in[0,t]$ on both sides of \r{e-mio} and estimating $\|v[\mu_\tau]-v[\nu_\tau]\|_{C^0}\leq K W_p(\mu_\tau,\nu_\tau)$, we have
\bqn
\phi(t)\leq e^{ 2L t}\phi(0) + \frac{Ke^{Lt}(e^{Lt}-1)}{L}\phi(t).
\eqnn
Take now $t'\leq \frac{1}{2L}$ and observe that $e^{Lt'}\leq 1+2Lt'\leq 2$. This gives
\bqn
\phi(t')\leq e\phi(0) + 4Kt\phi(t').
\eqnn
If moreover $t'\leq \frac{1}{8K}$, we have $\phi(t')\leq 2e \phi(0)$. Since such estimate holds for $t'\leq \min\Pg{\frac{1}{2L}\frac{1}{8K}}$, we apply it recursively $\Pceil{t \max\Pg{2L,8K}}$ times to find the estimate for a given $t$. Observing that $W_p(\mu_t,\nu_t)\leq \phi(t)$, we find \r{e-stimasol}.

Uniqueness of the solution is a direct application, recalling that $W_p(\mu,\nu)=0$ if and only if $\mu=\nu$.
\eproof

We finally prove an estimate about the convergence rate of \schema{Scheme 1} to the solution.
\bp
Let $v$ satisfy \Hp, and $\mu_0$ given. Let $\mu=\mu_{\Pq{0,T}}$ be the solution of \r{e-cauchy}, and $\mu^k$ the approximation of $\mu$ computed using \schema{Scheme 1} with $\dt=\frac{T}{2^k}$. Then, for a sufficiently small $\dt$, it holds
\bqn
d(\mu,\mu^k)\leq 4KM \dt\frac{e^{T (4L+4K)}-1}{4L+4K}.
\eqnl{e-schema1}
\ep
\bproof
The proof is similar to the proof of Proposition \ref{p-schema1}. For a time $t$, let $j$ be such that $t\in\Pq{j\dt,(j+1)\dt}$, and define $t'=t-j\dt\leq \dt$. We estimate $W_p(\mu_t,\mu^k_t)$ with respect to $W_p(\mu_{j\dt},\mu^k_{j\dt})$, independently on $p\geq 1$ and for sufficiently small $\dt$. Applying \r{e-stime4}, we have
\bqn
W_p(\mu_t,\mu^k_t)&=&W_p(\Phi^{v\Pq{\mu_t}}_{t'}\#\mu_{j\dt},\Phi^{v\Pq{\mu^k_{j\dt}}}_{t'}\#\mu^k_{j\dt})\leq\nn
&\leq& (1+4L\dt) W_p(\mu_{{j\dt}},\mu^k_{{j\dt}})+4\dt \sup_{\tau\in\Pq{{j\dt},t}}\|v\Pq{\mu_\tau}-v\Pq{\mu^k_{j\dt}}\|_{C^0}.
\eqnl{e-ora}
We estimate $\|v\Pq{\mu_\tau}-v\Pq{\mu^k_{{j\dt}}}\|_{C^0}$ using \Hp\ and \r{e-lipt}, as follows:
\bqn
\|v\Pq{\mu_\tau}-v\Pq{\mu^k_{{j\dt}}}\|_{C^0}&\leq& K \Pt{W_p(\mu_\tau,\mu_{{j\dt}})+W_p(\mu_{{j\dt}},\mu^k_{j\dt})}\leq \nn
&\leq& K \Pt{M\tau+W_p(\mu_{{j\dt}},\mu^k_{{j\dt}})}.
\eqnl{e-estv}
We to the supremum over $t\in\Pq{0,\dt}$ in \r{e-estv}. Replacing it in \r{e-ora}, we have
\bqn
W_p(\mu_t,\mu^k_t)&\leq& (1+C_1\dt)W_p(\mu_{{j\dt}},\mu^k_{{j\dt}})+\dt^2 C_2,
\eqnn
with $C_1:=4L+4K$ and $C_2:=4KM$. Applying it recursively in $j=0,\ldots, 2^k$, we have 
\bqn
d(\mu_t,\mu^k)\leq \dt^2 C_2 \frac{(1+C_1\dt)^{2^k}-1}{(1+C_1\dt)-1}\leq \dt C_2 \frac{e^{T C_1}-1}{C_1}.
\eqnn
\eproof

\section{Lagrangian and Eulerian numerical schemes}
\label{s-schemi}

In this section we consider other schemes for \r{e-cauchy}, with increasing level of discretization. We prove that all these schemes converge to the solution of \r{e-cauchy}, whose existence and uniqueness has just been proved. \schema{Scheme 2} is only an initial discretization in space, while \schema{scheme 3} is semidiscrete in time with initial discretization in space. \schema{Scheme 4} is Lagrangian discrete in space and time. We conclude definitions of schemes with \schema{scheme 5}, that is discrete Eulerian.

\schema{Schemes 2 and 3} can be considered as intermediate steps to the definition of \schema{schemes 4 and 5}, that are the main schemes practically used. Indeed, \schema{schemes 1, 2, 3} are not practically feasible in reality, since they need integration of ``general'' vector fields. Instead, \schema{schemes 4 and 5} can be used since they are based on approximations of such vector fields.

The last part of this section is devoted some simple generalizations of these results to equations used in modelling of pedestrian, namely several populations in interaction, eventually in a multi-scale approach.\\

We start this section by introducing spatial discretization. Fix a space discretization parameter $\h>0$ and discretize the space $\R^n$, dividing it in a grid of hypercubes with side $\h$. For simplicity, starting from now we will use the notation for grids in $\R^2$, calling each hypercube simply as ``square''. This also comes from our interest in pedestrian modelling, for which the space is usually 2-dimensional. We denote the grid with the symbol $\grid$. We also use the notation $\Chi_A$ for the characteristic function of the set $A\subset\R^n$. Given a measure $\mu$, we denote the discretized measure (with an abuse of notation) with the symbol $\grid\Pq{\mu}$, computed as follows: given a square $\square$ in the grid $\grid$, we define $\grid\Pq{\mu}$ on this square to be constant and with the same mass as $\mu$ on the square. More precisely, we have

\mmfunz{\grid}{\Pr}{\Pr}{\mu}{\grid\Pq{\mu}:=\displaystyle \sum_{\square\in\grid} \frac{\mu_0 (\square)}{\h^n} .}
This choice clearly permits to preserve the total measure, i.e. $\mu(\R^n)=\grid\Pq{\mu}(\R^n)$.

We now estimate $W_p(\mu,\grid\Pq{\mu})$, using a map $\gamma$ that redistributes the mass inside each square $\square$ only. Remark that $\gamma$ is not optimal in general. We have $W_p^p(\mu,\grid\Pq{\mu})\leq \int_{\R^n} |\gamma(x)-x|^p \,d\mu \leq \sqrt{n}^p \h^p$, i.e.
\bqn
W_p(\mu,\grid\Pq{\mu})\leq \sqrt{n} \h.
\eqnl{e-Wpgrid}
The same idea can be used to estimate $W_p(\mu,\grid\Pq{\mu})$ when $\mu$ is a positive measure with mass $\mu(\R^n)=C\neq 1$. We have 
\bqn
W_p(\mu,\grid\Pq{\mu})\leq \mu(\R^n)^{1/p} \sqrt{n} \h.
\eqnl{e-Wpgridmasse}

\subsection{Scheme with initial discretization in space}
\newcommand{\bmu}{\mu^{(2)}}

In this section we introduce a scheme with initial discretization in space. The scheme is based on a grid $\grid$ of squares with side $\h$. We use the notation $\bmu_t$ to denote the solution of this scheme. We define the initial data $\bmu_0$ as the discretization of the initial data $\mu_0$, i.e. $\bmu_0:=\grid\Pq{\mu_0}$. We then compute the approximated solution $\bmu_t$ given by evolution according to $v\Pq{\bmu_t}$. It is the exact solution starting from an approximated data. Existence and uniqueness for \r{e-cauchy} give us the good definition of $\bmu_t$. Remark that $\bmu_t$ is not constant over each square $\square$ of the initial grid, and moreover it is not even piecewise constant in general.

We now prove that this scheme is convergent.
\bp
Let $v$ satisfy \Hp, and $\mu_0$ given. Let $\bmu=\bmu_{\Pq{0,T}}$ be computed via \schema{Scheme 2} with parameter $\h$, and $\mu$ be the exact solution of \r{e-cauchy}. We have 
\bqn
d(\mu,\bmu)\leq \h \sqrt{n}\costantona{T}.
\eqnl{e-schema2}
Then, $\bmu$ converges weakly to $\mu$ in $C(\Pq{0,T},\Prr)$ for $\h\rightarrow 0$.
\ep
\bproof
Using \r{e-stimasol}, we have that $W_p(\bmu_t,\mu_t)\leq \costantona{t} W_p(\bmu_0,\mu_0)$, thus $$d(\bmu,\mu)\leq \costantona{T} W_p(\bmu_0,\mu_0).$$ Using \r{e-Wpgrid}, we have \r{e-schema2}. 
Since $d(\bmu,\mu)$ converges to $0$ for $\h\rightarrow 0$, then $\bmu$ converges weakly to $\mu$.
\eproof

\subsection{Semi-discrete in time Lagrangian scheme with initial discretization in space}

\renewcommand{\bmu}{\mu^{(3)}}
In this section we introduce a third numerical scheme, that is a semi-discrete in time Lagrangian scheme with initial discretization in space.

Fix a space discretization parameter $\h$ and a time discretization $\dt$. We perform a first approximation $\bmu_0$ of the initial data $\mu_0$ in space, that is $\bmu_0:=\grid\Pq{\mu_0}$. We then compute its evolution using \schema{Scheme 1}. The resulting scheme is the following.

\begin{pschema}{SCHEME 3\\Lagrangian, semi-discrete in time,\\ initial discretization in space}

\b{INITIALIZATION:} Fix a time discretization parameter $\dt$ and a space discretization parameter $\dx$.

Given the starting measure $\mu_0$, define $\bmu_0:=\grid\Pq{\mu_0}$.\\

%\b{STEP 0:} \\% Define $v_0$ on the interval $\Pq{0,\dt}$ as $v_0:=v\Pq{\bmu_0}$ and compute the corresponding flow $\flow^{v_0}_t$. 

%For $t\in\Pq{0,\dt}$, define $\bmu_t:=\flow^{v_0}_t\#\bmu_0$.\\

\b{STEP i+1:} Given $\bmu_{i\dt}$, define $v_{i\dt}$ on the interval $\Pq{i\dt,(i+1)\dt}$ as $v_{i\dt}:=v\Pq{\bmu_{i\dt}}$ and compute the corresponding flow $\flow^{v_{i\dt}}_t$. 

For $t\in\Pq{i\dt,(i+1)\dt}$, define $\bmu_t:=\flow^{v_{i\dt}}_t\#\bmu_{i\dt}$.\\

\b{STOP:} when reaching $T$. 
\end{pschema}

We now prove that also this scheme is convergent.
\bp
\label{p-schema3}
Let $v$ satisfy \Hp, and $\mu_0$ given. Let $\bmu=\bmu_{\Pq{0,T}}$ be computed via \schema{Scheme 3} with parameters $\h$ and $\dt=\frac{T}{2^k}$, and $\mu$ be the exact solution of \r{e-cauchy}. Then, for sufficiently small $\dt$, we have
\bqn
d(\mu,\bmu)\leq \h\sqrt{n}\costantona{T} +
4KM \dt\frac{e^{T (4L+4K)}-1}{4L+4K}.
\eqnn
In particular, $\bmu$ converges weakly to $\mu$ in $C(\Pq{0,T},\Prr)$ for $\h,\dt\rightarrow 0$.
\ep
\bproof
Define $\nu$ the exact solution of \r{e-cauchy} starting from $\bmu_0$, that is the approximated solution of $\mu$ given by \schema{Scheme 2} with parameter $h$. Thus, $\bmu$ is the approximated solution of $\nu$ given by \schema{Scheme 1} with parameter $\dt$.

Since $d(\mu,\bmu)\leq d(\mu,\nu) +d(\nu,\bmu)$, then using estimates \r{e-schema1}-\r{e-schema2}, we have the result.
\eproof

\subsection{Discrete Lagrangian scheme with velocity of centers}

\renewcommand{\bmu}{\mu^L}
In this section we introduce a fourth numerical scheme, that can be seen as an approximation of the previous one. We indicate with $\mu^L_t$ the function computed via this scheme. The index ${}^L$ means ``Lagrangian''. 

Fix a space discretization parameter $\h$ and a time discretization $\dt$. We perform a first approximation $\bmu_0$ of the initial data $\mu_0$ in space, that is $\bmu:=\grid\Pq{\mu_0}$. We take a square of the grid and denote it with $\square_0$. We denote the center of this square as $x_{\square_0}$ and $v(\square_0)$ to indicate the velocity $v\Pq{\bmu_0}(x_{\square_0})$, i.e. the vector field $v\Pq{\bmu_0}$ evaluated at the center of the square. This can be seen as the approximation of the vector field $v\Pq{\bmu_0}$ with a vector field that is piecewise-constant on the same grid.

We compute the evolution $\square_\dt$ of the square $\square_0$ as its translation of the vector $\dt\, v({\square_0})$. We define $\bmu_\dt$ as the sum of all the translated squares, i.e. 
$$\bmu_\dt:=\sum_{\square_0\in\grid} \bmu_0(\square_0) \frac{\Chi_{\square_\dt}}{\h^n}.$$ We repeat the translation $\square_{2\dt}$ of the square $\square_\dt$, following $v(\square_\dt):=v\Pq{\bmu_\dt}(x_{\square_\dt})$, i.e. the new velocity field evaluated at the center of the squares. We continue until reaching $T$.

\begin{pschema}{SCHEME 4\\Lagrangian, discrete in space and time,\\ velocity at the centers}

\b{INITIALIZATION:} Fix a space discretization parameter $\h$ and a time discretization parameter $\dt$.% Define a grid $\grid$ dividing the space $\R^n$ into squares of side $\h$.%\\

Given a starting measure $\mu_0$, define $\bmu_0$ as $\bmu_0:=\grid\Pq{\mu_0}$. \\
%For $t\in\Pq{0,\dt}$, define $\bmu_t:=\sum_{\square_0\in\grid} \bmu_0(\square_0) \frac{\Chi_{\square_t}}{\h^n}$ with $\square_t=\square_0+t\, v(\square_0)$.\\

\b{STEP i+1:} Given $\bmu_{i\dt}$, define for $t\in\Pq{i\dt,(i+1)\dt}$ $$\bmu_{t}:=\sum_{\square_0\in\grid} \bmu_{i\dt}(\square_0) \frac{\Chi_{\square_t}}{\h^n}$$
 with $\square_{t}=\square_{i\dt}+(t-i\dt)\, v(\square_{i\dt})$, where $v(\square_{i\dt}):=v\Pq{\bmu_{i\dt}}(x_{\square_\dt})$ and $x_{\square_{i\dt}}$ is the center of the square $\square_{i\dt}$.\\

\b{STOP:} when reaching $T$. 
\end{pschema}

\newcommand{\tmu}{\tilde{\mu}}
\newcommand{\mub}{\mu^{(3)}}
\newcommand{\tx}{\tilde{x}}
\newcommand{\supp}{{\mathrm{supp}}}

One of the advantages of \schema{Scheme 4} is that the evolution of a square is still a square, moreover of the same dimension and with sides parallel to axes of $\R^n$. In other words, the measure $\bmu_t$ is always piecewise-constant. Remark instead that the grid is used  at the beginning of the algorithm only, and that afterwards the squares do not belong to the grid. Moreover, the squares can overlap and follow different velocities, thus the evolution cannot be written globally as a solution of \r{e-cauchy}, since $v\Pq{\bmu_t}(x)$ is not uniquely determined.

Another important feature of the scheme is that the value of the function inside each square does not change. Nevertheless, there is an interaction among all the squares, since the velocity of each square is given by the map $v\Pq{\mu_t}$. One can also observe that, if at time $l\dt$ two different squares have the same center, then they have the same dynamics from that moment on. This is possible, since the evolution of centres is not given by a flow (that would give existence and uniqueness), but a discrete-time dynamics.

We now prove that the scheme is convergent. 
\bp
\label{p-schema4}
Let $v$ satisfy \Hp, and $\mu_0$ be given. Let $\bmu=\bmu_{\Pq{0,T}}$ be computed using \schema{Scheme 4} with parameters $\h$ and $\dt=\frac{T}{2^k}$. Let $\mu$ be the exact solution of \r{e-cauchy}. Then $d(\mu,\bmu)\leq d(\mu,\mu^{(3)}) + d(\mu^{(3)},\bmu)$ and, for $\dt<\frac{\log(2)}{L}$, it holds
\bqn
d(\mu^{(3)},\bmu)\leq 2L \sqrt{n}\frac{e^{KT(e^{LT})}-1}{K} \h.
\eqnn

Then $\bmu$ converges weakly to $\mu$ in $C(\Pq{0,T},\Prr)$ for $\h,\dt\rightarrow 0$.
\ep
\bproof The idea of the proof is to decompose the evolution of measures $\mub_t$ and $\bmu_t$ in the evolution of each square of the initial grid.

Take a square $\square_0$ at the starting grid and consider the measure $\tmu_0:= \mub_0 \Chi_{\square_0}$.  Let $t\in \Pq{j\dt,(j+1)\dt}$ and call $S:=j\dt$, $t'=t-S$. Consider the two following evolutions: the first is the evolution given by $\tmu_t:=\Phi^{v\Pq{\mub_S}}_{t'}\# \tmu_S$. The second is the evolution given by $\nu^L_t:=\Phi^{v(\square_S)}_{t'}\# \nu^L_S$. Observe that in both formulas, the evolution is computed according to the vector field given by the whole measures $\mub_S$ and $\bmu_S$. As a consequence, we have $\mub_t=\sum_{\square_t} \tmu_t$ and $\mu^L_t=\sum_{\square_t} \nu^L_t$.

We denote with $\tilde{x}_t$ the evolution of the center of the square $\square_t$ according to the first evolution, that is $\tilde{x}_t=\Phi^{v\Pq{\mub_S}}_{t'} (\tilde{x}_S)$. The evolution of the center of the square $\square_t$ according to the second evolution is $x^L_t=\Phi^{v\Pq{\bmu_S}(x^L_S)}_{t'}(x^L_S)$.

We estimate the distance $W_p(\tmu_t, \nu^L_t)$, using the decomposition represented in Figure \ref{fig-schema4}.\\

\newcommand{\storto}[5]{
\pgfmoveto{\pgfxy(#1)}
\pgfcurveto{\pgfxy(#2)}{\pgfxy(#2)}{\pgfxy(#3)}
\pgfmoveto{\pgfxy(#3)}
\pgfcurveto{\pgfxy(#4)}{\pgfxy(#4)}{\pgfxy(#5)}
\pgfstroke}

\begin{figure}[htb]
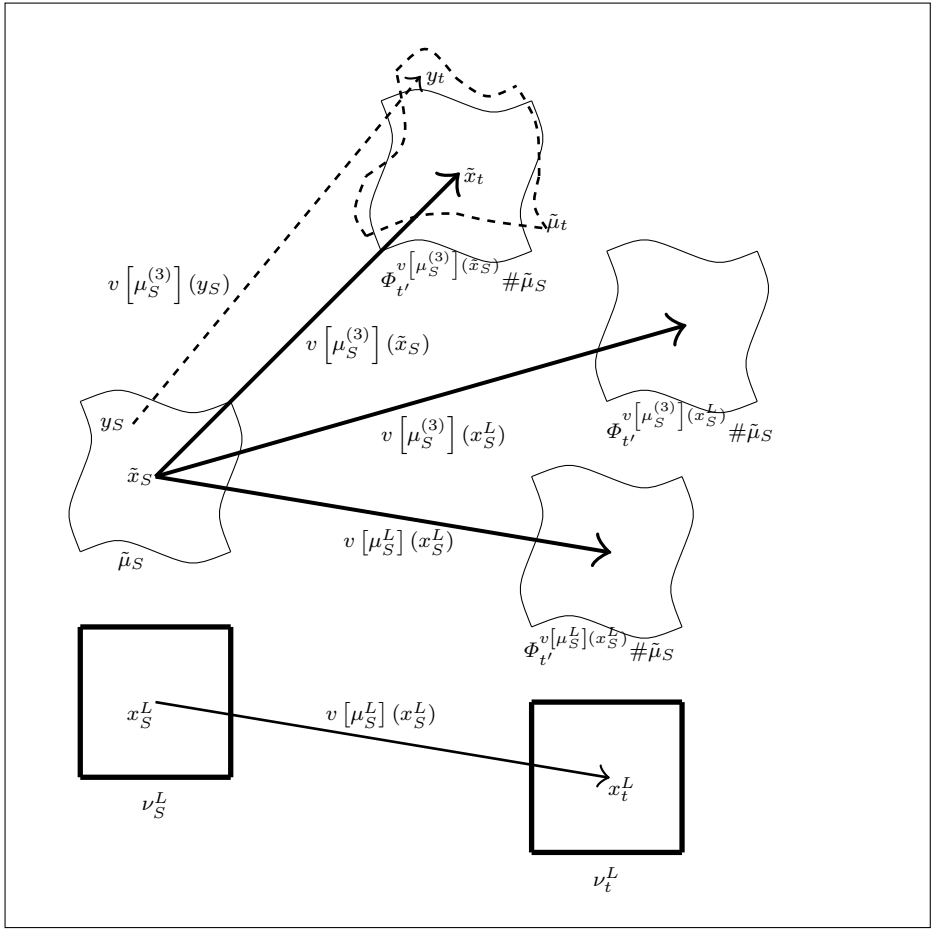


\begin{center}
\begin{pgfpictureboxed}{0cm}{0cm}{12.3 cm}{12.3cm}
%\pgfgrid{\pgfxy(0,0)}{\pgfxy(12,12)}

\pgfputat{\pgfxy(1.5,5)}{\pgfbox[left,top]{$\tmu_S$}}
\pgfputat{\pgfxy(1.8,6)}{\pgfbox[center,center]{$\tx_S$}}
\storto{1,5}{.8,5.5}{1,6}{1.2,6.5}{1,7}
\storto{1,5}{1.5,5.2}{2,5}{2.5,4.8}{3,5}
\storto{3,5}{2.8,5.5}{3,6}{3.2,6.5}{3,7}
\storto{1,7}{1.5,7.2}{2,7}{2.5,6.8}{3,7}

\pgfputat{\pgfxy(6.9,4)}{\pgfbox[left,top]{$\Phi^{v\Pq{\mu^L_S}(x^L_S)}_{t'}\#\tmu_S$}}
\storto{7,4}{6.8,4.5}{7,5}{7.2,5.5}{7,6}
\storto{7,4}{7.5,4.2}{8,4}{8.5,3.8}{9,4}
\storto{9,4}{8.8,4.5}{9,5}{9.2,5.5}{9,6}
\storto{7,6}{7.5,6.2}{8,6}{8.5,5.8}{9,6}

\pgfputat{\pgfxy(8,7)}{\pgfbox[left,top]{$\Phi^{v\Pq{\mub_S}(x^L_S)}_{t'}\#\tmu_S$}}
\storto{8,7}{7.8,7.5}{8,8}{8.2,8.5}{8,9}
\storto{8,7}{8.5,7.2}{9,7}{9.5,6.8}{10,7}
\storto{10,7}{9.8,7.5}{10,8}{10.2,8.5}{10,9}
\storto{8,9}{8.5,9.2}{9,9}{9.5,8.8}{10,9}

\pgfputat{\pgfxy(8-3,9)}{\pgfbox[left,top]{$\Phi^{v\Pq{\mub_S}(\tx_S)}_{t'}\#\tmu_S$}}
\storto{8-3,7+2}{7.8-3,7.5+2}{8-3,8+2}{8.2-3,8.5+2}{8-3,9+2}
\storto{8-3,7+2}{8.5-3,7.2+2}{9-3,7+2}{9.5-3,6.8+2}{10-3,7+2}
\storto{10-3,7+2}{9.8-3,7.5+2}{10-3,8+2}{10.2-3,8.5+2}{10-3,9+2}
\storto{8-3,9+2}{8.5-3,9.2+2}{9-3,9+2}{9.5-3,8.8+2}{10-3,9+2}

\begin{pgfscope}
\pgfsetlinewidth{1pt}
\pgfsetdash{{3pt}{3pt}}{0pt}
\storto{7.8-3,7.2+2}{7.6-3,7.5+2}{7.8-3,8+2}{8.4-3,8.4+2}{8.2-3,9.4+2}
\storto{8.2-3,9.4+2}{8.6-3,9.8+2}{9-3,11.3}{9.4-3,9+2}{9.8-3,9.2+2}
\storto{7.8-3,7.2+2}{8.5-3,7.5+2}{9-3,7.5+2}{9.4-3,7.4+2}{10.2-3,7.3+2}
\storto{10.2-3,7.3+2}{10-3,7.6+2}{10.1-3,8+2}{10.1-3,8.7+2}{9.8-3,9.2+2}
\pgfputat{\pgfxy(7.2,9.5)}{\pgfbox[left,top]{$\tmu_t$}}
\pgfputat{\pgfxy(6.1,10)}{\pgfbox[left,center]{$\tx_t$}}
\begin{pgfscope}
\pgfsetendarrow{\pgfarrowto}
\pgfline{\pgfxy(1.7,6.7)}{\pgfxy(5.5,11.3)}
\pgfputat{\pgfxy(3,8.3)}{\pgfbox[right,bottom]{$v\Pq{\mub_S}(y_S)$}}
\pgfputat{\pgfxy(1.6,6.7)}{\pgfbox[right,center]{$y_S$}}
\pgfputat{\pgfxy(5.6,11.3)}{\pgfbox[left,center]{$y_t$}}
\end{pgfscope}
\end{pgfscope}

\begin{pgfscope}
\pgfsetendarrow{\pgfarrowto}
\pgfsetlinewidth{1.5pt}
\pgfline{\pgfxy(2,6)}{\pgfxy(8,5)}
\pgfputat{\pgfxy(4.5,5.3)}{\pgfbox[left,top]{$v\Pq{\mu^L_S}(x^L_S)$}}
\pgfline{\pgfxy(2,6)}{\pgfxy(8+1,5+3)}
\pgfputat{\pgfxy(5,6.8)}{\pgfbox[left,top]{$v\Pq{\mub_S}(x^L_S)$}}
\pgfline{\pgfxy(2,6)}{\pgfxy(6,10)}
\pgfputat{\pgfxy(4,8)}{\pgfbox[left,top]{$v\Pq{\mub_S}(\tx_S)$}}
\end{pgfscope}

\begin{pgfscope}
\pgfsetlinewidth{2pt}
\pgfputat{\pgfxy(2,1.8)}{\pgfbox[center,top]{$\nu^L_S$}}
\pgfputat{\pgfxy(1.8,3)}{\pgfbox[center,top]{$x^L_S$}}
\pgfline{\pgfxy(1,2)}{\pgfxy(1,4)}
\pgfline{\pgfxy(1,2)}{\pgfxy(3,2)}
\pgfline{\pgfxy(3,2)}{\pgfxy(3,4)}
\pgfline{\pgfxy(3,4)}{\pgfxy(1,4)}

\pgfputat{\pgfxy(8,.8)}{\pgfbox[center,top]{$\nu^L_t$}}
\pgfputat{\pgfxy(8.2,2)}{\pgfbox[center,top]{$x^L_t$}}
\pgfline{\pgfxy(7,1)}{\pgfxy(7,3)}
\pgfline{\pgfxy(7,1)}{\pgfxy(9,1)}
\pgfline{\pgfxy(9,1)}{\pgfxy(9,3)}
\pgfline{\pgfxy(9,3)}{\pgfxy(7,3)}

\begin{pgfscope}
\pgfsetendarrow{\pgfarrowto}
\pgfsetlinewidth{1pt}
%\pgfsetdash{{3pt}{3pt}}{0pt}
\pgfline{\pgfxy(2,3)}{\pgfxy(8,2)}
\pgfputat{\pgfxy(5,3)}{\pgfbox[center,top]{$v\Pq{\mu^L_S}(x^L_S)$}}
\end{pgfscope}
\end{pgfscope}

\end{pgfpictureboxed}

\caption{Convergence of \schema{Scheme 4}: decomposition.}
\label{fig-schema4}
\end{center}

\end{figure}

We have
\bqn
W_p(\tmu_t, \nu^L_t)&=&
W_p(\Phi^{v\Pq{\mub_S}}_{t'}\#\tmu_S,\Phi^{v\Pq{\mu^L_S}(x^L_S)}_{t'}\#\nu^L_S)
\leq
\eqnn
\bqn
&\leq&
W_p(\Phi^{v\Pq{\mub_S}}_{t'}\#\tmu_S,\Phi^{v\Pq{\mub_S}(\tx_S)}_{t'}\#\tmu_S)+
W_p(\Phi^{v\Pq{\mub_S}(\tx_S)}_{t'}\#\tmu_S,\Phi^{v\Pq{\mub_S}(x^L_S)}_{t'}\#\tmu_S)+\nn
&+&
W_p(\Phi^{v\Pq{\mub_S}(x^L_S)}_{t'}\#\tmu_S,\Phi^{v\Pq{\mu^L_S}(x^L_S)}_{t'}\#\tmu_S)+
W_p(\Phi^{v\Pq{\mu^L_S}(x^L_S)}_{t'}\#\tmu_S,\Phi^{v\Pq{\mu^L_S}(x^L_S)}_{t'}\#\nu^L_S).
\eqnn
We now estimate the four terms.

We estimate the first term using \r{e-stime4massa}. We estimate $\|v\Pq{\mub_S}-v\Pq{\mub_S}(\tx_S)\|_{C^0}$ on the support of $\tmu_S$ only. Take a point $y_S\in\supp\Pt{\tmu_S}$ and estimate $|v\Pq{\mub_S}(y_S)-v\Pq{\mub_S}(\tx_S)|$. Since $y_S\in\supp\Pt{\tmu_S}$, then it is the evolution of a point $y_0\in\supp\Pt{\tmu_S}$ via \schema{Scheme 3}, thus
\bqn 
|y_{(j+1)\dt}-\tx_{(j+1)\dt}|&=&|y_{j\dt}+\dt\, v\Pq{\mub_{j\dt}}(y_{j\dt})-\Pt{\tx_{(j+1)\dt}+\dt\, v\Pq{\mub_{j\dt}}(\tx_{j\dt})}|\leq\nn
&\leq& |y_{j\dt}-\tx_{j\dt}|+\dt L|y_{j\dt}-\tx_{j\dt}|,
\eqnn
hence $|y_S-\tx_S|=|y_{j\dt}-\tx_{j\dt}|\leq |y_0-\tx_0| (1+\dt L)^j$. Since $y_0$ belongs to the square $\square_0$ and $\tx_0$ is its center, we have $|y_0-\tx_0|\leq \frac{\sqrt{n}}{2}\h$. Recall that $v\Pq{\mub_S}$ is Lipschitz, thus $\|v\Pq{\mub_S}-v\Pq{\mub_S}(\tx_S)\|_{C^0}\leq L \frac{\sqrt{n}}{2}\h (1+\dt L)^j$. Summing up, we have
\bqn
W_p(\Phi^{v\Pq{\mub_S}}_{t'}\#\tmu_S,\Phi^{v\Pq{\mub_S}(\tx_S)}_{t'}\#\tmu_S)\leq 
\tmu_S(\R^n)^{1/p} e^{L t'/p}(e^{Lt'}-1)\frac{\sqrt{n}}{2}\h (1+\dt L)^j.
\eqnn

We estimate the second and third terms. Remark that, since a vector field evaluated at one point is a vector, then each of these two terms represents the Wasserstein distance of the same measure under two different translations. It is easy to prove that, given two vectors $a,b$, it holds 
$$W_p(\Phi^a_t\# \mu,\Phi^b_t\#\mu)\leq \Pt{\int |(x+t\,a)-(x+t\,b)|^p \,d\mu}^{1/p}=|t||a-b|\mu(\R^n)^{1/p}.$$
This is indeed a particular case of \r{e-stime4massa}. For the second term, we have
\bqn
W_p(\Phi^{v\Pq{\mub_S}(\tx_S)}_{t'}\#\tmu_S,\Phi^{v\Pq{\mub_S}(x^L_S)}_{t'}\#\tmu_S)
&\leq&
t' L |\tx_S-x^L_S | \tmu_S(\R^n)^{1/p}.
\eqnn
For the third term, recalling the Lipschitzianity of $v$, we have 
\bqn
W_p(\Phi^{v\Pq{\mub_S}(x^L_S)}_{t'}\#\tmu_S,\Phi^{v\Pq{\mu^L_S}(x^L_S)}_{t'}\#\tmu_S)
&\leq&
t' K W_p(\mub_S,\bmu_S) \tmu_S(\R^n)^{1/p}.
\eqnn

The fourth term has a particular structure too, since it is the Wasserstein distance of two measures under the same translation $t'\,v\Pq{\mu^L_S}(x^L_S)$. Since the Wasserstein distance is invariant under translation, we have
\bqn
W_p(\Phi^{v\Pq{\mu^L_S}(x^L_S)}_{t'}\#\tmu_S,\Phi^{v\Pq{\mu^L_S}(x^L_S)}_{t'}\#\nu^L_S)=
W_p(\tmu_S,\nu^L_S).
\eqnn

Summing up the four estimates, we have
\bqn
&&W_p(\tmu_t, \nu^L_t)\leq
\tmu_S(\R^n)^{1/p} e^{L t'/p}(e^{Lt'}-1)\frac{\sqrt{n}}{2}\h (1+\dt L)^{j}+\label{e-dopo4}\\
&&\hspace{3mm}+
t' L |\tx_S-x^L_S | \tmu_S(\R^n)^{1/p}+
t' K W_p(\mub_S,\bmu_S) \tmu_S(\R^n)^{1/p}+
W_p(\tmu_S,\nu^L_S).
\eqnn
The estimate is increasing as a function of $t'$, thus the supremum is attained for $t'=\dt$. We now estimate $|\tx_{(j+1)\dt}-x^L_{(j+1)\dt} |$ in terms of $W_p(\mub_{j\dt},\bmu_{j\dt})$. Recall that $\tx_t$ and $x^L_t$ are the evolution of the same starting point $x^L_0$ under two different numerical schemes. Thus
\bqn
|\tx_{(j+1)\dt}&-&x^L_{(j+1)\dt} |=|\tx_{j\dt}+\dt\,v\Pq{\mub_{j\dt}}(\tx_{j\dt})-(x^L_{j\dt}+\dt\,v\Pq{\bmu_{j\dt}}(x^L_{j\dt}))|\leq\nn
&&\leq |\tx_{j\dt}-x^L_{j\dt}| +\dt \left(|v\Pq{\mub_{j\dt}}(\tx_{j\dt})-v\Pq{\mub_{j\dt}}(x^L_{j\dt})|+\right.\nn
&&+ \left. |v\Pq{\mub_{j\dt}}(x^L_{j\dt})-v\Pq{\bmu_{j\dt}}(x^L_{j\dt})|\right)\leq\nn
&&\leq |\tx_{j\dt}-x^L_{j\dt}|+L\dt |\tx_{j\dt}-x^L_{j\dt}| + K\dt W_p(\mub_{j\dt},\bmu_{j\dt}).
\eqnn
One can observe that $|\tx_{0}-x^L_{0}|=0$, and also $W_p(\mub_{0},\bmu_{0})=0$, since the two schemes have same initial data. Thus $|\tx_{\dt}-x^L_{\dt}|=0$. We thus have 
\bqn
|\tx_{j\dt}-x^L_{j\dt} | &\leq& \frac{(1+L\dt)^j-1}{(1+L\dt)-1} K\dt\sup_{i\in\Pg{0,1,\ldots, j-1}} W_p(\mub_{i\dt},\bmu_{i\dt}) \leq\nn
&\leq & \frac{K}{L}(e^{jL\dt}-1)\sup_{i\in\Pg{0,1,\ldots,j-1 }} W_p(\mub_{i\dt},\bmu_{i\dt}).
\eqnn
We plug this formula in \r{e-dopo4}.

We now pass to the global estimate of the distance. Fix a time $t$ and for each square $\square$ with $\mu_0(\square)>0$, compute the following quantity: $W_p(\tmu_t, \nu^L_t)\tmu_S(\R^n)^{-1/p}$. Now maximise this quantity among all squares and denote it with $\phi(t)$. Now come back to \r{e-dopo4}. We multiply for $\tmu_S(\R^n)^{-1/p}$ on both sides, then maximise on all squares on both sides. It gives
\bqn
&&\phi((j+1)\dt)\leq
 e^{L \dt/p}(e^{L\dt}-1)\frac{\sqrt{n}}{2}\h (1+\dt L)^{j}+\nn
&&%\hspace{-3cm}
+ \dt K (e^{jL\dt}-1) \mbox{$\sup_{i\in\Pg{0,1,\ldots,j-1 }}$} W_p(\mub_{i\dt},\bmu_{i\dt}) +
\dt K W_p(\mub_S,\bmu_S) +\phi(j\dt).
\eqnn
Now call $f_j:=\sup_{\tau\in\Pq{0,j\dt}} \phi(t)$. Since $W_p(\mub_t,\bmu_t)\leq \phi(t)$ due to \r{e-Wpsup}, we have 
\bqn
f_{j+1}\leq e^{L \dt/p}(e^{L\dt}-1)\frac{\sqrt{n}}{2}\h (1+\dt L)^{j}+
(1+\dt K e^{jL\dt}) f_j.
\eqnn
We maximize the coefficients $(1+\dt L)^{j}$ and $e^{jL\dt}$ with respect to $j\leq 2^k$, that give $e^{LT}$ in both cases. We estimate the term $f_{2^k}$ under the hypothesis $L\dt< \log(2)$, for which we have $e^{L\dt}< 1+ 2L\dt$ and $e^{L\dt/p}\leq e^{L\dt}<2$. Since $f_0=0$, we have
\bqn
f_{2^k}\leq e^{L \dt/p}(e^{L\dt}-1)\frac{\sqrt{n}}{2}\h e^{LT} \frac{(1+\dt K e^{LT})^{2^k}-1}{1+\dt K e^{LT}-1}\leq
2L \sqrt{n}\h\frac{e^{KT(e^{LT})}-1}{K}.
\eqnn
Since we have $d(\mub,\bmu)\leq f_{2^k}$, we have the estimate. We thus have weak convergence of $\bmu$ to $\mub$ for $\h\rightarrow 0$. Using Proposition \ref{p-schema3} under the additional hypothesis $\dt\rightarrow 0$, we have weak convergence of $\mub$ to $\mu$, then of $\bmu$ to $\mu$.
\eproof

\brem
The convergence of the Lagrangian schemes is not really surprising in the context of pedestrian modelling by using measures. Indeed, we approximate a set of discrete agents (the pedestrians) with a continuous measure. This is the passage from microscopic to macroscopic model. Then, the scheme is convergent for $\h\rightarrow 0$, that means that the we go back to a microscopic scale, i.e. the pedestrians.
\erem

One can improve this scheme by computing a more precise evolution of each square. For example, one can allow deformations of the axes, rotations and so on. This idea coincide with the idea of computing the vector field $v$ for each square as a certain approximation of the original $v$. For our scheme, we simply evaluate $v$ at the center, i.e. we perform a Taylor expansion of order $0$. Improvements of this kind certainly result in better convergence rates, but need a more complicated implementation.

%Indeed, even this scheme has a drawback, since at each step it is necessary to sum all the squares to find $\bmu_t$ and compute $v\Pq{\bmu_t}$. One of the goals of the following scheme is exactly to solve this drawback.

\subsection{Eulerian scheme}
\renewcommand{\bmu}{\mu^E}
\renewcommand{\tmu}{\tilde{\mu}}
\newcommand{\tnu}{\tilde{\nu}}

In this section we present a last scheme to compute numerically the solutions of \r{e-cauchy}. This scheme has been first proposed in \cite{ben-ARMA} for modelling of pedestrians. We indicate the solution of the scheme as $\mu^E_t$. The index ${}^E$ means ``Eulerian''. We call it Eulerian, since we are interested in the evolution of the value $\mu_t(P)$ at a point $P$ not changing in time. Instead, the previous schemes where Lagrangian, since we were interested on the spatial evolution of a point $x$, i.e. its trajectory $x(t)$. This was particularly clear in the two previous schemes, in which we fixed a starting square $\square_0$ and study its evolution in time $\square_t$.

Fix a space discretization parameter $\h>0$ and a time discretization $\dt$. We perform a first approximation $\bmu_0$ of the initial data $\mu_0$ in space, that is $\bmu_0:=\grid\Pq{\mu_0}$. Take a square of the grid $\square\in\grid$. We denote the center of this square as $x_{\square}$ and $v_t(\square)$ to indicate the velocity $v\Pq{\bmu_t}(x_{\square})$, i.e. the vector field $v\Pq{\bmu_t}$ evaluated at the center of the square.

We compute the evolution of the square $\square$ as its translation of the vector $\dt\, v_0({\square})$, i.e. $\square+\dt\, v_0({\square})$. We define $\tmu_\dt$ as the sum of all the translated squares, i.e. $\tmu_\dt:=\sum_{\square\in\grid} \bmu_0(\square) \frac{\Chi_{\square+\dt\, v_0({\square})}}{\h^n}$. For the moment, this coincide with \schema{Scheme 4}. The difference is that we define $\bmu_\dt$ as the approximation of $\tmu_\dt$ computing the mean values on the starting grid, i.e. $\bmu_\dt=\grid\Pq{\tmu_\dt}$.

We then repeat the same idea ``evolution + mean values'' starting from $\bmu_\dt$. We continue until reaching $T$.

\begin{pschema}{SCHEME 5\\Eulerian, discrete in space and time,\\ velocity at the centers}

\b{INITIALIZATION:} Fix a space discretization parameter $\h$ and a time discretization parameter $\dt$.% Define a grid $\grid$ dividing the space $\R^n$ into squares of side $\h$.%\\

Given a starting measure $\mu_0$, define $\bmu_0$ as $\bmu_0:=\grid\Pq{\mu_0}$. \\
%For $t\in\Pq{0,\dt}$, define $\bmu_t:=\sum_{\square_0\in\grid} \bmu_0(\square_0) \frac{\Chi_{\square_t}}{\h^n}$ with $\square_t=\square_0+t\, v(\square_0)$.\\

\b{STEP i+1:} Define $$\tmu_{(i+1)\dt}:=\sum_{\square\in\grid} \bmu_{i\dt}(\square) \frac{\Chi_{\square+ \dt\, v_{i\dt}(\square)}}{\h^n}$$ with $v_{i\dt}(\square):=v\Pq{\bmu_{i\dt}}(x_{\square})$ and $x_{\square}$ is the center of the square $\square$.

Define the approximated solution $\bmu_{(i+1)\dt}:=\grid\Pq{\tmu_{(i+1)\dt}}$.\\

\b{STOP:} when reaching $T$.
\end{pschema}

We now prove that the scheme is convergent. First observe that $\mu^E$ is defined for times $j\dt$ only, thus the distance $d(\mu,\mu^E)$ is not defined by equation \r{e-d}. Nevertheless, we redefine it with a slight abuse of notation $$d(\mu,\mu^E):=\sup\Pg{W_p(\mu_t,\mu^E_t) \mbox{ for } t\in\Pq{0,T}\mbox{ such that } \mu^E_t \mbox{ is defined.}}.$$

We have the following result. 
\bp
\label{p-schema5}
Let $v$ satisfy \Hp, and $\mu_0$ be given. Let $\bmu$ be computed using \schema{Scheme 5} with parameters $\h$ and $\dt=\frac{T}{2^k}$. Let $\mu$ be the exact solution of \r{e-cauchy}. Then $d(\mu,\bmu)\leq d(\mu,\mu^L) + d(\mu^L,\bmu)$ and we have

\bi
\i for $p>1$:
\bqn
d(\mu^L,\mu^E)\leq \frac{(2^{1-1/p} \dt L +1)e^{T K}\sqrt{n}}{2^{1-1/p} -1} \h 2^{(1-1/p)\frac{T}{\dt}}.
\eqnl{e-schema5a}
For $\dt\rightarrow 0$ and $\lim_{\h,\dt} \h 2^{(1-1/p)\frac{T}{\dt}}=0$, $\mu^E$ converges weakly to $\mu$;

\i for $p=1$:
\bqn
d(\mu^L,\mu^E)\leq \frac{\sqrt{n}( \dt L +1) (e^{TK}-1)}{ K} \frac{\h}{\dt}.
\eqnl{e-schema5b}
For $\dt\rightarrow 0$ and $\lim_{\h,\dt} \frac{\h}{\dt}=0$, $\mu^E$ converges weakly to $\mu$.
\ei
\ep
\bproof
We want to estimate $W_p(\mu^L_{j\dt},\bmu_{j\dt})$. For this reason, we rewrite the algorithm in a different way. Assume to have already computed $\mu^L$ and $\mu^E$ on the whole interval $\Pq{0,T}$. Now fix a square $\square_0$ from the starting grid and follow the evolution of the mass inside $\square_0$ according to the two schemes. We denote the two evolutions at time $j\dt$ with $\nu^L_j$ and $\nu^E_j$ respectively. Remark that the support of $\nu^L_j$ is always a square: we denote such square with $\square^L_j$ and $x^L_j$ its center. The support of $\nu^E_j$ is instead the union of a finite number $m_j$ of squares: we denote such squares with $\square^E_{j,i}$, indexed by $i=1,\ldots, m_j$. They all belong to the grid $\grid$. We denote with $x^E_{j,i}$ the center of each square. For technical reasons, we always assume that $\square^L_j\subset \cup_{i=1}^{m_j} \square^E_{j,i}$, eventually adding squares with no mass to the set of squares $\square^E_{j,i}$.

We have $\nu^L_0=\nu^E_0=\frac{\mu_0(\square_0)}{\h^n}\chi_{\square_0}$, and 
\bqn
\nu^L_{j+1}&=&\Phi^{v\Pq{\mu^L_{j\dt}}(x^L_j)}_\dt \#\nu^L_j,\qquad
\nu^E_{j+1}=\grid\Pt{\tnu_{j+1}},
\eqnn
with
\bqn
\tnu_{j+1}&=&\sum_{i=1}^{m_j} \tnu_{j+1,i}\mbox{~~~and~~~}
\tnu_{j+1,i}:=\Phi^{v\Pq{\mu^E_{j\dt}}(x^E_{j,i})}_\dt \#(\nu^E_j \chi_{\square^E_{j,i}}).
\eqnn

One can rewrite the evolutions as follows. Define the vector fields 
\bqn
v^a(x):=
\begin{cases}
v\Pq{\mu^L_{j\dt}}(x^L_j) &\mbox{~~if~~} x\in \square^L_j,\\
v\Pq{\mu^E_{j\dt}}(x^E_{j,i}) &\mbox{~~if~~} x\in \square^E_{j,i}\backslash \square^L_j,\\
0 &\mbox{otherwise}
\end{cases}\nn
v^b(x):=
\begin{cases}
v\Pq{\mu^E_{j\dt}}(x^E_{j,i}) &\mbox{~~if~~} x\in \square^E_{j,i},\\
0 &\mbox{otherwise}.
\end{cases}
\eqnn
We can rewrite $\nu^L_{j+1}=\Phi^{v^a(x^L_j)}_\dt\#\nu^L_j$, $\tnu_{j+1,i}=\Phi^{v^b(x^E_{j,i})}_\dt\#\nu^E_{j,i}$  and $\nu^E_{j+1}=\grid\Pt{\tnu_{j+1}}$. Remark that, as already stated, none of the two evolutions is given by a global flow, thus we cannot directly apply results of Section \ref{s-Wpflow}. Nevertheless, take $\gamma$ to be the optimal map realizing $W_p(\nu^L_j,\nu^E_j)$ with $\gamma\#\nu^E_j=\nu^L_j$ and define $\nu^L_{j,i}:=\gamma\#\nu^E_{j,i}$. Observe that 
\bqn
W^p_p(\nu^L_{j+1},\Phi^{v^b}_\dt\#\nu^E_j)\leq \sum_{i=1}^{m_j} W^p_p(\Phi^{v^a(x^L_j)}_\dt\#\nu^L_{j,i},\Phi^{v^b(x^E_{j,i})}_\dt\#\nu^E_{j,i}).
\eqnl{e-add}
We estimate the right hand side, that is the Wassertein distance under translations of constant vectors $v^a(x^L_j)$ and $v^b(x^E_{j,i})$ respectively. Remark that $\gamma$ is the optimal map realizing $W_p(\nu^L_{j,i},\nu^E_{j,i})$, by construction. Then it holds 
\bqn
&&W^p_p(\Phi^{v^a(x^L_j)}_\dt\#\nu^L_{j,i},\Phi^{v^b(x^E_{j,i})}_\dt\#\nu^E_{j,i})\leq\nn
&&\leq\int_{\R^n} |x - \Pt{\gamma(x-\dt v^b(x^E_{j,i}))+\dt v^a(x^L_j)})|^p\,d\Pt{\Phi^{v^b(x^E_{j,i})}_\dt\#\nu^E_{j,i}}(x)\leq\nn
&&\leq\int_{\R^n} |y +\dt v^b(x^E_{j,i})- \Pt{\gamma(y)+\dt v^a(x^L_j)})|^p\,d\nu^E_{j,i}(y)\leq\nn
&&\leq 2^{p-1} \Pt{ \int_{\R^n} |y -\gamma(y)|^p \,d\nu^E_{j,i}(y) +\dt^p|v^b(x^E_{j,i})- v^a(x^L_j)|^p \nu^E_{j,i}(\R^n)}=\nn
&&= 2^{p-1} \Pt{W_p^p(\nu^L_{j,i},\nu^E_{j,i})+ \dt^p|v^b(x^E_{j,i})- v^a(x^L_j)|^p \nu^E_{j,i}(\R^n)}.
\eqnl{e-adesso}
Observe that 
\bqn
| v^b(x^E_{j,i})- v^a(x^L_j)|&\leq&\sup \Pg{|v\Pq{\mu^E_{j\dt}}(x^E_{j,i}) - v\Pq{\mu^L_{j\dt}}(x^L_j)|\ \mbox{s.t.}\ \square^E_{j,i}\cap \square^L_j\neq \emptyset},
\eqnn
where we have used that $\square^L_j\subset \cup_{i=1}^{m_j} \square^E_{j,i}$. Remark that $\square^E_{j,i}\cap \square^L_j\neq \emptyset$ when $|x^L_j-x^E_{j,i}|\leq \sqrt{n} \h$. Hence
\bqn
| v^b(x^E_{j,i})- v^a(x^L_j)|&\leq&\sup \Pg{\scriptstyle |v\Pq{\mu^L_{j\dt}}(x^L_j)- v\Pq{\mu^E_{j\dt}}(x^L_j)|+ |v\Pq{\mu^E_{j\dt}}(x^L_j)-v\Pq{\mu^E_{j\dt}}(x^E_{j,i})\ |\ \square^E_{j,i}\cap \square^L_j\neq \emptyset}\leq\nn
&\leq&  K W_p(\mu^L_{j\dt},\mu^E_{j\dt})+L \sqrt{n} \h.
\eqnn
Plug this estimate in \r{e-adesso} and sum over all $i$ in \r{e-add}. Also recall that, by contruction, we have $W_p^p(\nu^L_{j},\nu^E_{j})=\sum_{i=1}^{m_j} W_p^p(\nu^L_{j,i},\nu^E_{j,i})$. Then compute the $p$-th root, recalling that $(a+b)^{1/p}\leq a^{1/p}+b^{1/p}$. It gives
\bqn
W^p(\nu^L_{j+1},\Phi^{v^b}_\dt\#\nu^E_j)&\leq& 2^{1-1/p} (W_p(\nu^L_{j},\nu^E_{j})+ \dt K W_p(\mu^L_{j\dt},\mu^E_{j\dt})\nu^E_{j}(\R^n)^{1/p}+\nn
&&~~~~+\dt L \sqrt{n} \h \nu^E_{j}(\R^n)^{1/p}).
\eqnn
Now observe that 
\bqn
&&W_p(\nu^L_{j+1},\nu^E_{j+1})\leq W_p(\nu^L_{j+1},\tnu_{j+1})+ W_p(\tnu_{j+1},\nu^E_{j+1})\leq\nn
&&\leq 2^{1-1/p} \Pt{W_p(\nu^L_{j},\nu^E_{j})+ \dt K W_p(\mu^L_{j\dt},\mu^E_{j\dt})\nu^E_{j}(\R^n)^{1/p}+\dt L \sqrt{n} \h \nu^E_{j}(\R^n)^{1/p}}+\nn
&&+ \sqrt{n} \h \nu^E_{j}(\R^n)^{1/p},
\eqnl{e-questa}
where we have estimated $W_p(\tnu_{j+1},\nu^E_{j+1})=W_p(\tnu_{j+1},\grid\Pt{\tnu_{j+1}})$ via \r{e-Wpgridmasse}. Multiply on the left side for $\nu^E_{j+1}\Pt{\R^n}^{-1/p}$ and on the right side for $\nu^E_{j}\Pt{\R^n}^{-1/p}$, that coincide.
Define $\phi_j:=\sup_{\square_0\in\grid} W_p(\nu^L_j,\nu^E_j) \nu^L_j\Pt{\R^n}^{-1/p}$. Then \r{e-questa} reads as
\bqn
\phi_{j+1}\leq 2^{1-1/p} \phi_j+ 2^{1-1/p} \dt K \phi_j +2^{1-1/p} \dt L\sqrt{n} \h +\sqrt{n} \h,
\eqnn
where we have used $W_p(\mu^L_{j\dt},\mu^E_{j\dt})\leq \phi_j$, as proved in \r{e-Wpsup}. Recall that $\phi_0=0$ and observe that $2^{1-1/p} + 2^{1-1/p} \dt K > 1$. Then, for $p>1$, we have
\bqn
d(\mu^L,\nu^E)&\leq& \sup_{j=0,1,\ldots,2^k}\phi_{j}\leq \sqrt{n} \h (2^{1-1/p} \dt L +1) \frac{(2^{1-1/p} + 2^{1-1/p} \dt K)^{2^k}-1}{2^{1-1/p} + 2^{1-1/p} \dt K-1}\leq\nn
&\leq& \sqrt{n} \h (2^{1-1/p} \dt L +1) \frac{2^{(1-1/p)\frac{T}{\dt}} e^{2^k \dt K}}{2^{1-1/p} -1}\leq\nn
&\leq& \frac{(2^{1-1/p} \dt L +1)e^{T K}\sqrt{n}}{2^{1-1/p} -1} \h 2^{(1-1/p)\frac{T}{\dt}}.
\eqnl{e-finale}
This gives estimate \r{e-schema5a} and convergence  to $\mu^L$ for $\lim_{\h,\dt} \h 2^{(1-1/p)\frac{T}{\dt}}=0$.

The case $p=1$ gives 
\bqn
d(\mu^L,\nu^E)&\leq& \sqrt{n} \h ( \dt L +1) \frac{(1 + \dt K)^{2^k}-1}{1 + \dt K-1}\leq \sqrt{n}( \dt L +1) \frac{e^{TK}-1}{\dt K}\h.
\eqnn
This gives estimate \r{e-schema5b} and the convergence to $\mu^L$ for $\lim_{\h,\dt} \frac{\dx}{\dt}=0$.

\noindent The convergence of $\mu^E$ to $\mu$ is then given by Proposition \ref{p-schema4}, since $\mu^L$ converges to $\mu$.
\eproof

\brem The condition given above can be surprising if compared with the well-known CFL condition for convergence of schemes for the transport equation. We recall that CFL asks for $\lim_{\dx,\dt} \frac{\dx}{\dt}>0$, see e.g. \cite{ref-CFL}. A comparison between these two conditions is outside the goals of this paper, nevertheless we give two remarks concerning this issue.

First, it is important to state that Proposition \ref{p-schema5} is a partial result, since we can't find a counterexample of non-convergence of \schema{Scheme 5} under hypothesis \Hp.

On the other hand, CFL is used for schemes in which the distance between measures is the $L^1$-distance. In our setting this is not a natural choice, for reasons that are explained in detail in the following Section \ref{s-L1}. In other words, a scheme that is convergent with respect to one distance can be non-convergent with respect to another one.
\erem

\brem
It is interesting to observe that, if $\lim_{\h,\dt} \h 2^{(1-1/p)\frac{T}{\dt}}=0$ is satisfied for a certain $p_1$, then it is satisfied for all $p_2<p_1$. Moreover, $\lim_{\h,\dt} \frac{\dx}{\dt}=0$, i.e. the condition for $p_2=1$ is satisfied too. This is a direct consequence of the fact that condition \Hp\ becomes stronger for decreasing $p$, thus the corresponding condition for convergence of the scheme can be weaker.
\erem

%a point $x\in\R^n$ and consider the displacement of this point under \schema{Scheme 4}. Its displacement is given by $v\Pq{\mu}(\tx)$ with $\tx$ the center of the square to which $x$ belongs. Since $|\tx-x|\leq \frac{\sqrt{2}}{2}h$, we have $|v\Pq{\mu}(\tx)-v\Pq{\mu}(x)|\leq Lh$. Given a square $
%\bqn
%W_p(\mu^L_{(l+1)\dt},\tmu_{(l+1) \dt})&\leq&\displaystyle \sum_{\square_0\in\grid}W_p(\Phi^{v\Pq{\mu^L_{l\dt}}(\tx^_S)}_\dt\#\mu^L_{l\dt}(\square_0) \Chi_{\square_{l\dt}},
%W_p(T\#\mu^L_{l\dt}(\square_0) \Chi_{\square_{l\dt}}.
%\eqnn

\subsection{Generalizations}
\label{s1-generalizations}

In this section we present two simple generalizations of the previous results of convergence of numerical schemes. We first focus on the generalization of the problem \r{e-cauchy} to a problem involving several populations. We then focus on problems in which the populations are no more absolutely continuous with respect to the Lebesgue measure, but contain Dirac delta too. This second generalization is used for multi-scale modelling, see e.g. \cite{ben-multiscale}.

\subsubsection{Generalization to several populations}
\newcommand{\muu}{{\boldsymbol \mu}}
\newcommand{\nuu}{{\boldsymbol \nu}}

In this section we show how to pass from the problem \r{e-cauchy} for a population $\mu$ to the case of a finite number of populations $\mu^1,\ldots,\mu^N$. We explain in the following the precise mathematical setting that we are going to study. We already observe that this generalization is fundamental for the applications to pedestrian modelling, in which one very often  has  interaction among populations with different goals (e.g. pedestrians interacting at a pedestrian crossing).

We first precisely define the problem. Given $N\in\N$ and $T>0$, consider $\muu:=\Pt{\mu^1,\ldots,\mu^N}\in C(\Pq{0,T},(\Prr)^N)$ each of them satisfying
\bqn
\begin{cases}
\partial_t\mu^i+\nabla\cdot(v^i\Pq{\muu} \mu^i)=0,\\
\mu^i_{|_{t=0}}=\mu^i_0.
\end{cases}
\label{e-cauchyN}
\eqn
In other words, each of the components $\mu^i$ satisfies its continuity equation as in \r{e-cauchy}, but the vector field $v^i\Pq{\muu}$ depends on the global $\muu$.

In this context, it is necessary to define distances both for vectors and for measures in $\Pr^N$. Given $v,w\in (R^n)^N$, we define $|v-w|:=\sum_{i=1}^N |v^i-w^i|$. Similarly, given $\muu,\nuu\in\Pr^N$, we define $$W_p(\muu,\nuu):=\sum_{i=1}^N W_p(\mu^i,\nu^i).$$ Remark that this is one of the possible choices of the distance, corresponding to use 1-norm in $\R^N$. All the following results hold true by changing the norm.

We also have to change the hypotheses $\Hp$ to the following
\newcommand{\HN}{{\bf (H-N)}}

\vspace{3mm}

\noindent\begin{Sbox}\begin{minipage}{\textwidth}\vspace{2mm}\begin{center}\HN\vspace{5mm}\\
\begin{minipage}{0.9\textwidth}
The function
\mmfunz{v\Pq{\muu}}{\Pr^N}{\Pt{C^{1}(\R^n)\cap L^\infty(\R^n)}^N}{\muu}{v\Pq{\muu}}
satisfies 
\bi
\i $v\Pq{\muu}$ is uniformly Lipschitz and uniformly bounded, i.e. there exist $L$, $M$ not depending on $\muu$, such that for all $\muu\in(\Pr)^N, x,y\in\R^n,$
\bqn
\hspace{-5mm}|v\Pq{\muu}(x)-v\Pq{\muu}(y)|\leq L |x-y|\qquad |v\Pq{\muu}(x)|\leq M.
\eqnn

\i $v$ is a Lipschitz function, i.e. there exists $K$ such that 
\bqn\|v\Pq{\muu}-v\Pq{\nuu}\|_{\mathrm{C^0}} \leq K \Wp{\muu,\nuu}.
\eqnn
\ei
\end{minipage}\end{center}\vspace{2mm}\end{minipage}\end{Sbox}\fbox{\TheSbox}\vspace{3mm}

We now review the main results of the paper given for $N=1$ under hypothesis $\Hp$. We show that the same results hold for general $N$, under hypothesis $\HN$. We start with the convergence of \schema{Scheme 1}. The definition of this scheme (and of the other ones presented in this paper) is a straightforward adaptation of the scheme on each component.

\bp
\label{p-schema1N}
Let $v$ satisfy \HN, and $\muu_0\in(\Prr)^N$ be given. Let $\muu^k=\muu_{\Pq{0,T}}^k$ be constructed by \schema{Scheme 1} with $\dt=\frac{T}{2^k}$. Then, the sequence $\Pg{\muu_{\Pq{0,T}}^k}_{k\in\N}$ is a Cauchy sequence in $C(\Pq{0,T},(\Pr)^N)$.
\ep
\bproof
The proof is completely equivalent to the proof of Proposition \ref{p-schema1}. Using the same notations, one has to estimate $\|v^k_j-v^{k+1}_{2j}\|_{C^0}\leq K W_p(\muu^k_j,\muu^{k+1}_{2j})$. Rewrite estimates for each component. Then, summing over all components and rearranging the terms, one finds the same recursive formula
$$W_p(\mu^k_t,\mu^{k+1}_t)\leq (1+2^{-k} C_1) W_p(m^k_j, m^{k+1}_{2j})+2^{-2k}C_2,$$
where $C_1$ and $C_2$ coincide with the ones given in the proof of Proposition \ref{p-schema1}, after replacing $K$ by $K\,N$. This gives convergence of the series.
\eproof

Since $C(\Pq{0,T},(\Pr)^N)$ is complete, then it exists $\muu=\lim_k \muu^k$. This limit is a solution of \r{e-cauchyN} and belongs to $C(\Pq{0,T},(\Prr)^N)$. This can be proved by working on each component as in Theorem \ref{t-esistenza}. We have 
\bt
\label{t-esistenzaN}
Let $v$ satisfy \HN, and $\muu_0\in(\Prr)^N$ be given. Let $\muu^k=\muu^k_{\Pq{0,T}}\in C(\Pq{0,T},(\Prr)^N)$ be computed using \schema{Scheme 1} with $\dt=\frac{T}{2^k}$, starting with $\muu^k_0=\muu_0$. Then the limit $\bar{\muu}=\lim_k \muu^k$ exists and is a solution of \r{e-cauchyN}. Moreover, this solution is unique.
\et

The uniqueness comes from the same estimate \r{e-stimasol}, in which constants have to be computed more precisely. Similarly, all results of convergence of schemes can be adapted to this case, computing the new constants.

\subsubsection{Generalization to multi-scale processes}
\label{s-multiscala}

In this section we show how to pass from the problem \r{e-cauchyN} for $N$ populations $\mu^i$, each of them absolutely continuous, to populations presenting Dirac delta. We do not study the general problem in which a Dirac delta can be transformed in an absolutely continuous measure. Rather, we consider evolutions in which the absolutely continuous and the singular parts cannot melt. In this sense, each measure $\mu_t$ can be always be written as 
\bqn
\mu_t=\mu^{ac}_t +\mu^d_t=\gamma^{ac}_t\#\mu^{ac}_0+\gamma^{d}_t\#\mu^d_0,
\eqnl{e-acd}
where $\mu^{ac}_t$ is the absolutely continuous part, and $\mu^d_t$ is the sum of the Dirac delta.

This strong condition permits to study this problem as the problem of two distinct measures $\mu^{ac}$ and $\mu^d$, each of them satisfying its own transport equation \r{e-cauchyN}, with $v\Pq{\mu}=v\Pq{\mu^{ac},\mu^d}$. We recall that the transport equation for a Dirac delta $\delta_x$ is equivalent to the ODE $\dot x= v(x)$. Thus one can use \schema{Scheme 1} to find an approximated solution of this problem, then Proposition \ref{p-schema1N} to check convergence, then verify that the limit is indeed a solution. The unique detail to verify is that estimates used in the proof of Proposition \ref{p-schema1N} hold for this special case with $\mu_t=\gamma^{ac}_t\#\mu^{ac}_0+\gamma^{d}_t\#\mu^d_0$. Given this decomposition, one can estimate the Wasserstein distance $W_p(\mu,\nu)$ by estimating $W_p(\mu^{ac},\nu^{ac})$ and $W_p(\mu^d,\nu^d)$. Thus, the only result to prove is that estimates given in Section \ref{s-Wpflow} also hold for measures that are sum of Dirac delta only. A careful look to the proofs show that they are verified, mainly using the Gronwall lemma for the ODE $\dot x= v(x)$.

Finally, it is evident that the same results hold in the case of $N$ populations, each of them satisfying the decomposition rule \r{e-acd}.

\section{Definitions of $v$ for modelling of pedestrian flow}
\label{s-nuclei}

In this section we present some choices for the definition of $v=v\Pq{\mu}$ that have been proposed for the modelling of pedestrian flows. We briefly describe the general idea, referring to \cite{ben-ARMA} for more details.

Take a pedestrian inside a crowd of people with spatial density $\mu$. He has a goal (reach a place, get out of a building, etc...), that defines a vector field $v^d(x)$, that is the field of {\it desired velocity}. Trajectories to reach the goal are thus integral curves of this field. If the pedestrian is alone, his dynamics is only determined by $v=v^d$. If instead the pedestrian is part of a crowd, it is good for him to locally deviate from the trajectory given by $v^d$ to avoid crowded areas (that would force him to reduce his velocity). Hence, one can add a {\it  velocity of interaction} $v^i$ that represents the tendency of avoiding crowded areas.  It is clear that $v^i$ depends on the actual $\mu$, i.e. $v^i=v^i\Pq{\mu}$.

Thus, a model for the evolution of crowds is given by \r{e-cauchy} with $$v=v\Pq{\mu}=v^d+v^i\Pq{\mu}.$$ We want now study the regularity of $v$, according to some models that have been proposed for $v^i$. As a first observation, the regularity depends on the regularity of both terms $v^d$ and $v^i$. We always assume that $v^d$ is Lipschitz and bounded. Under this hypothesis, $v$ satisfies \Hp\ if and only if $v^i$ does. Remark that the problem of regularity of $v^d$ is non-trivial, for example in presence of domains with obstacles, see e.g. \cite{maury}.

From now on, we can consider $v=v^i$. Take a pedestrian in $x\in\R^n$ and consider how he avoids crowded areas. One can imagine that he checks the density of the crowd in a finite area around him and stays away of crowded areas. A simple model is the following: take a kernel $\fz{\eta}{\R^n}{[ 0,+\infty )}$, representing the weight of interaction of the single pedestrian with people in his \nei. For example, it is reasonable to model that overcrowded places close to the pedestrian are ``more important'' or better detected than far ones. Define 
$$x^*:=\frac{\int_{\R^n}y\,\eta(x-y)\,d\mu(y)}{\int_{\R^n}\eta(x-y)\,d\mu(y)},$$ that is indeed the center of mass of the crowd with respect to the kernel $\eta$, and define
\bqn
v\Pq{\mu}(x):=(x-x^*) f\Pt{\int_{\R^n}\eta(x-y)\,d\mu(y)}.
\eqnl{e-v1}
where $f$ is a non-decreasing weight function. This formula represent the fact that $v$ drives away from $x^*$, with an intensity given by $f$. The meaning of $f$ non-descreasing is the following: if the area is crowded, i.e. $\int_{\R^n}\eta(x-y)\,d\mu(y)$ is big, then the pedestrian gets away of the crowd faster than if the area is not crowded. Observe that \r{e-v1} is not defined for $\int_{\R^n}\eta(x-y)\,d\mu(y)=0$. In this case, we simply define $v\Pq{\mu}(x)=0$. This choice does note change the dynamics, since an area without mass undergoes no change in time.

We now study two particular cases of velocity \r{e-v1}. The first is given by $f(x)\equiv 1$, while the second is given by $f(x)=x^\al$ with $\al\geq 1$. We prove that in the first case $v$ does not satisfy \Hp, while in the second it does.

\begin{figure}[htb]
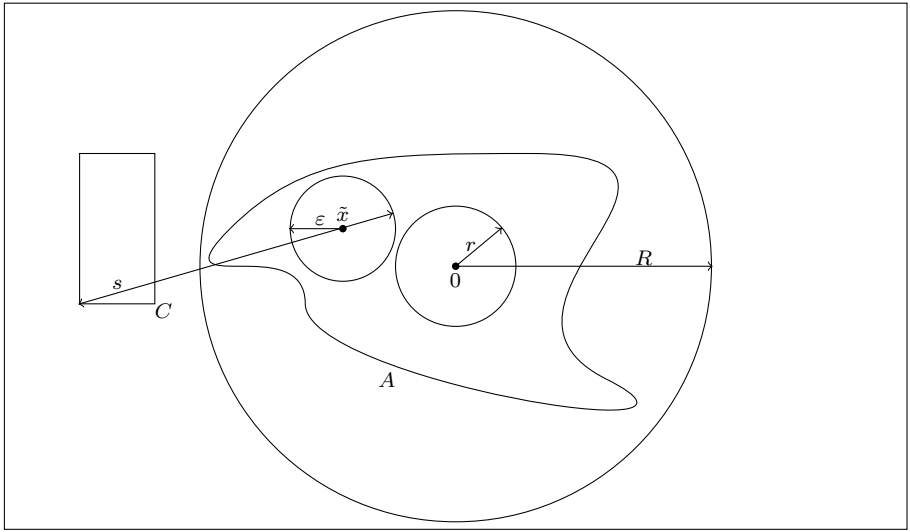

\begin{center} \begin{pgfpictureboxed}{0cm}{0cm}{12 cm}{7 cm}
\newcommand{\alto}{+1}
%\pgfgrid[stepx=1cm,stepy=1cm]{\pgfxy(0,0)}{\pgfxy(12,6)}
\pgfxycurve(4,2\alto)(4,3\alto)(2,2\alto)(3,3\alto)
\pgfxycurve(3,3\alto)(4,4\alto)(5,4\alto)(7,4\alto)
\pgfxycurve(7,4\alto)(10,4\alto)(6,2\alto)(8,1\alto)
\pgfxycurve(8,1\alto)(10,0\alto)(4,1\alto)(4,2\alto)

\pgfcircle[stroke]{\pgfxy(6,2.5\alto)}{.8cm}
\pgfcircle[fill]{\pgfxy(6,2.5\alto)}{.05cm}
\pgfputat{\pgfxy(6,3.4)}{\pgfbox[center,top]{$0$}}
\pgfcircle[stroke]{\pgfxy(6,2.5\alto)}{3.4cm}

\pgfcircle[fill]{\pgfxy(4.5,4)}{.05cm}
\pgfputat{\pgfxy(4.5,4.1)}{\pgfbox[center,bottom]{$\tilde x$}}
\pgfcircle[stroke]{\pgfxy(4.5,4)}{.7cm}

\pgfrect[stroke]{\pgfxy(1,3)}{\pgfxy(1,2)}
\pgfputat{\pgfxy(2,3)}{\pgfbox[left,top]{$C$}}

\begin{pgfscope}
\pgfsetendarrow{\pgfarrowto}
\pgfline{\pgfxy(6,3.5)}{\pgfxy(9.4,3.5)}
\pgfputat{\pgfxy(8.5,3.52)}{\pgfbox[center,bottom]{$R$}}
\pgfline{\pgfxy(6,3.5)}{\pgfxy(6.6,4)}
\pgfputat{\pgfxy(6.2,3.7)}{\pgfbox[center,bottom]{$r$}}
\pgfline{\pgfxy(4.5,4)}{\pgfxy(3.8,4)}
\pgfputat{\pgfxy(4.2,4.05)}{\pgfbox[center,bottom]{$\eps$}}
\end{pgfscope}

\begin{pgfscope}
\pgfsetstartarrow{\pgfarrowto}
\pgfsetendarrow{\pgfarrowto}
\pgfline{\pgfxy(1,3)}{\pgfxy(5.15,4.2)}
\pgfputat{\pgfxy(1.5,3.2)}{\pgfbox[center,bottom]{$s$}}
\end{pgfscope}

\pgfputat{\pgfxy(5.2,1.9)}{\pgfbox[right,bottom]{$A$}}

\end{pgfpictureboxed}

\caption{The convolution kernel does not satisfy \Hp\ for $f\equiv 1$.}
\label{fig-contronucleo}
\end{center}

\end{figure}

Let $f(x)\equiv 1$. Given a kernel $\eta$ that is not identically 0, we find a family of measures $\mu_t$ such that $v\Pq{\mu}$ is not even a continuous vector field with respect to $\mu_t$ in the Wasserstein space. The idea is explained in Figure \ref{fig-contronucleo}.\\

 Let $R$ be such that $\mathrm{supp}\Pt{\eta}\subset B_R(0)$. Since $\eta$ is continuous, we have that the set $A:=\Pg{\eta>0}$ is open. It is always possible to choose $r$ sufficiently small to have $A\backslash \overline{B}_r(0)$ nonempty. Since it is open, we can always choose a point $\tilde x$ in this set, and $\eps>0$ sufficiently small to have $B_\eps(\tilde x)\subset A\backslash B_r(0)$. Clearly, if $y\in B_\eps(\tilde x)$, then $|y|\geq r>0$. Finally, define a compact set $C$ of non-zero Lebesgue measure outside the ball $B_R(0)$, and $s$ the maximum distance between elements of $B_\eps(\tilde x)$ and $C$, i.e. $s=\sup \Pg{|x-y|\,\mbox{~ s.t.~}\, x\in B_\eps(\tilde x),\,y\in C}$. It is clear that $s<\infty$.

We are now ready to define the family $\mu_t$ of measures. Define $$\mu_t:=\Pt{t \frac{\Chi_{B_\eps(\tilde x)}}{\lam(B_\eps(\tilde x))} + (1-t) \frac{\Chi_{C}}{\lam(C)}}\,\lam,$$ where $\lam$ is the Lebesgue measure. Observe that $v\Pq{\mu_0}(0)=0$, since $\int_{\R^n} \eta(-y)\, d\mu_0(y)=\frac{1}{\lam(C)} \int_C 0\, d\lam(y)=0$. Now observe that, for $t>0$, we have $\int_{\R^n} \eta(-y), d\mu_t(y)>0$, hence
\bqn
|v\Pq{\mu_t}(0)|&=&\Pabs{\frac{\int_{\R^n} y\eta(-y), d\mu_t(y)}{\int_{\R^n} \eta(-y), d\mu_t(y)}}=\Pabs{
\frac{ \frac{t}{\lam(B_\eps(\tilde x))}
\int_{B_\eps(\tilde x)} y\eta(-y)\, d\lam(y)}{
\frac{t}{\lam(B_\eps(\tilde x))}
\int_{B_\eps(\tilde x)} \eta(-y)\, d\lam(y)}}\geq\nn
&\geq& \frac{\inf \Pg{|y|\, \mbox{~s.t.~}y\in B_\eps(\tilde x)} \int_{B_\eps(\tilde x)} \eta(-y)\, d\lam(y)}{\int_{B_\eps(\tilde x)} \eta(-y)\, d\lam(y)}\geq r.
\eqnn
As a consequence, $v\Pq{\mu_t}(0)$ is not continuous with respect to the parameter $t$. We now show that $\mu_t$ is continuous with respect to parameter $t$ in 0, i.e. $\lim_{t}W_p(\mu_0,\mu_t)=0$. Fix a time $t$ and consider the measure $\nu_t$ shared by $\mu_0$ and $\mu_t$, that is exactly $\nu_t:=(1-t) \frac{\Chi_{C}}{\lam(C)}\,\lam$. Thus, by \r{e-Wpshared}, we have
\bqn
W_p(\mu_0,\mu_t)=W_p(\mu_0-\nu_t,\mu_t-\nu_t)=W_p \Pt{t \frac{\Chi_{C}}{\lam(C)}\lam,t \frac{\Chi_{B_\eps(\tilde x)}}{\lam(B_\eps(\tilde x))}\lam}.
\eqnn
Take the optimal map $\gamma$ between these measures, and observe that $|x-\gamma(x)|\leq s$. Thus $W_p(\mu_0,\mu_t)\leq s t ^{1/p}$. Since $s$ and $\lam(C)$ are finite and do not depend on $t$, then $W_p(\mu_0,\mu_t)$ is continuous at $t=0$. Thus, $v\Pq{\mu_t}(0)$ is not continuous with respect to the distance $W_p(\mu_0,\mu_t)$.

We now study the case $f(x)=x^\al$ with $\al\geq 1$.
\bp
Let $v=v\Pq{\mu}$ defined by \r{e-v1}, with $\eta$ a positive, Lipschitz, bounded function with bounded support. Let $f(x)=x^\al$ with $\al\geq 1$. Then $v$ satisfies \Hp.
\ep
\bproof
Let $L$ be the Lipschitz constant of $\eta$, $M=|\eta|_\infty$ its maximal value, $R$ the radius of its bounded support, i.e. $\mathrm{supp}(\eta)\subset B_R(0)$. Call $\phi(x):=\int_{\R^n}\eta(x-y)\,d\mu(y)$. A direct computation shows that $|v\Pq{\mu}(x)|\leq R M^\al$. Similarly, we have
\bqn
|v\Pq{\mu}(x)-v\Pq{\mu}(z)| &=& \left|\int_{\R^n}(x-y)\eta(x-y)^\al\,d\mu(y)- \int_{\R^n}(z-y)\eta(z-y)^\al\,d\mu(y)\right|\leq\nn
&\leq& M^{\al-1} R \int_{\R^n} |\eta(x-y)-\eta(z-y)|\,d\mu(y)\leq M^{\al-1} R L |x-z|.
\eqnn
We finally prove that $v$ is Lipschitz with respect to the $W_1$ distance. We have
\bqn
|v\Pq{\mu}(x)-v\Pq{\nu}(x)|&\leq& M^{\al-1} \Pabs{\int_{\R^n}(x-y)\eta(x-y)\,d(\mu-\nu)(y)}.
\eqnn
Observe that $f(y):=(x-y)\eta(x-y)$ is a Lipschitz function, since
\bqn
|(x-y_1)\eta(x-y_1)-(x-y_2)\eta(x-y_2)|\leq R |\eta(x-y_1)-\eta(x-y_2)|\leq R L |y_1-y_2|.
\eqnn
Thus $\int_{\R^n}(x-y)\eta(x-y)\,d(\mu-\nu)(y)\leq R L W_1(\mu,\nu)$ via the Kantorovich-Rubinstein duality formula, hence $\|v\Pq{\mu}-v\Pq{\nu}\|_{C^0}\leq M^{\al-1} R L W_1(\mu,\nu)$.
\eproof

\section{Wasserstein vs $L^1$ distance}
\label{s-L1}

In this section we discuss the choice of the Wasserstein distance $W_p$ as the distance in the space $\Prr$. Indeed, this space is endowed with several other distances, see e.g. \cite{rachev}. A distance which plays an important role in this context is the $L^1$ distance, that is
\bqn
\|\mu-\nu\|_{L^1}:= \int |\mu(x)-\nu(x)| dx.
\eqnn
We will show that the choice of the Wasserstein distance is better than the choice of $L^1$ distance, both for modeling  and theoretical reasons.

For the modelling, observe that the Wasserstein metric is more adapted than $L^1$ distance to measure if two pedestrian populations are close or far. Indeed, take 3 different measures as in Figure \ref{fig:wassL1}, and call $\mu_i$ the measure centered in $x_i$, for $i=0,1,2$.

\immagine{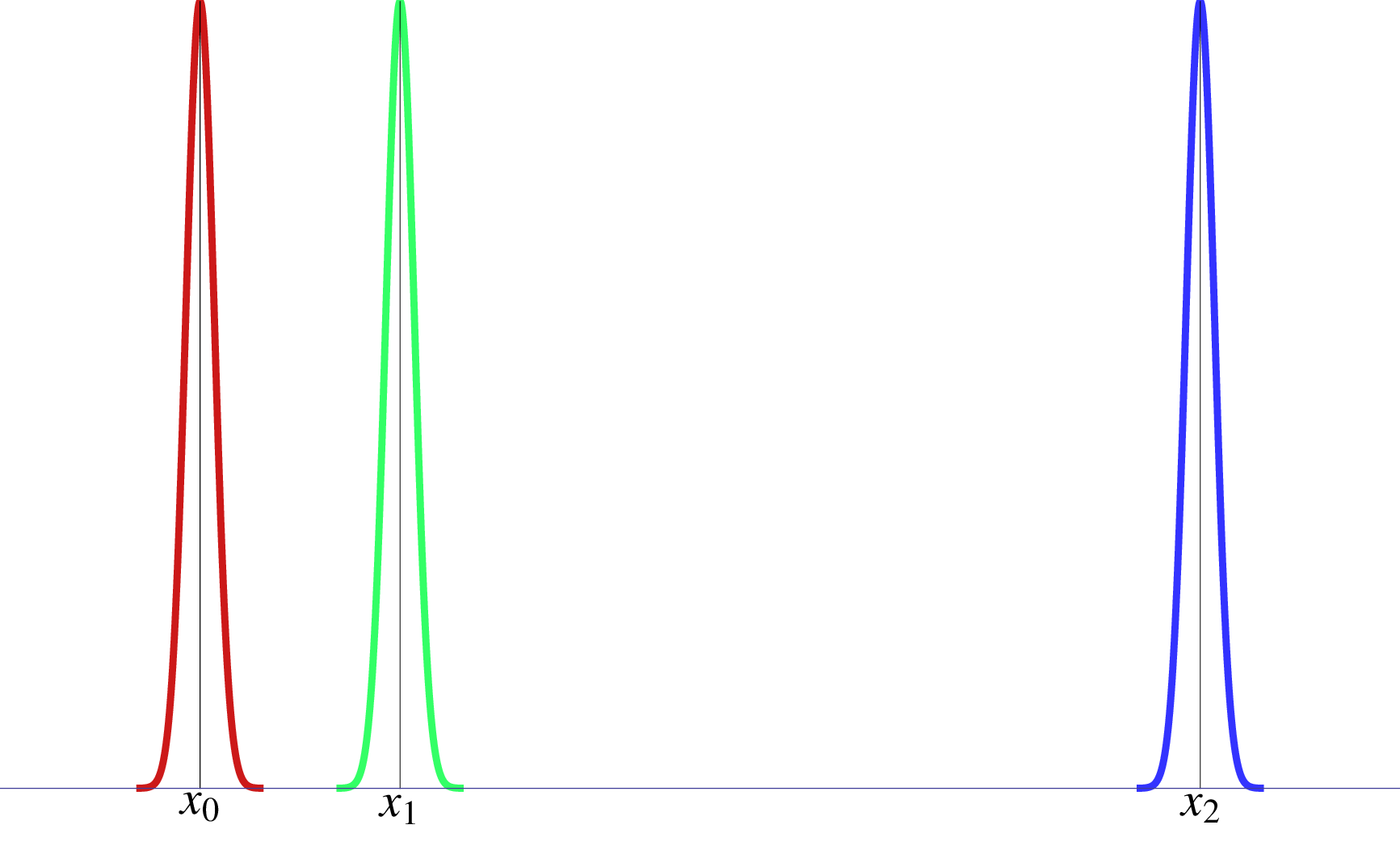}{Definition of $\mu_0,\mu_1,\mu_2$.}

It is easy to prove that
\bqn
W_p(\mu_i,\mu_j)=|x_i-x_j|,\qquad \mbox{~~~~and~~~~}\|\mu_i-\mu_j\|_{L^1}=2 \mbox{~~for~~} i\neq j.
\eqnn
In particular, the $L^1$ distance between the measures $\mu_0$ and $\mu_1$ is identical to the distance between $\mu_0$ and $\mu_2$. This is not natural with respect to our perception of distance between pedestrian crowds, that is better modeled by Wasserstein distance.

Another interesting issue comes from the kind of velocities usually used in the context of modelling of pedestrians flows. As showed in Section \ref{s-nuclei}, we often deal with velocities defined by convolution with a certain kernel. In this case, two measures that are close with respect to the $W_p$ distance but far with respect to $L^1$ distance give velocity fields that are close. For example, define $\mu_\eps:=\frac{1}{2\eps} \chi_{\Pt{-\eps,\eps}}$ and observe that $\mu_\eps$ and $\mu_{2\eps}$ are close for small $\eps$ with respect to distance $W_p$, but not with respect to $L^1$ distance.

On the other hand, if one has a velocity $v=v\Pq{\mu}$ defined by the value of $\mu$ in a point (provided a good definition of this quantity), i.e. $v\Pq{\mu}(x):=f(\mu(x))$ for a certain $f$, then the example given above provides a completely different result. Choose $f(\mu(x))=\mu(x)$. We have $v\Pq{\mu_\eps}(0)=2 v\Pq{\mu_{2\eps}}(0)$. Thus velocities are far, even if the measures are close with respect to the $W_p$ distance.

From the mathematical point of view, we show that, in our context, Lipschitzianity with respect to $L^1$ distance is not a good condition for  the problem \r{e-cauchy} with $\mu_0\in\Prr$. In particular, we show that it does not guarantee uniqueness of the solution. Indeed, assume the following hypotheses on $v$:
\newcommand{\Hpp}{{\bf ($\mathbf H_1$)}}

\vspace{3mm}
\noindent\begin{Sbox}\begin{minipage}{\textwidth}\vspace{2mm}\begin{center}\Hpp\vspace{5mm}\\
\begin{minipage}{0.9\textwidth}
The function
\mmfunz{v\Pq{\mu}}{\Pr}{C^{1}(\R^n)\cap L^\infty(\R^n)}{\mu}{v\Pq{\mu}}
satisfies 
\bi
\i $v\Pq{\mu}$ is uniformly Lipschitz and uniformly bounded, i.e. there exist $L$, $M$ not depending on $\mu$, such that for all $\mu\in\Pr, x,y\in\R^n,$
\bqn
\hspace{-5mm}|v\Pq{\mu}(x)-v\Pq{\mu}(y)|\leq L |x-y|\qquad |v\Pq{\mu}(x)|\leq M.
\eqnn

\i $v$ is a $L_1$-Lipschitz function, i.e. there exists $K$ such that 
\bqn\|v\Pq{\mu}-v\Pq{\nu}\|_{\mathrm{C^0}} \leq K \|\mu-\nu\|_{L^1}.
\eqnn
\ei
\end{minipage}\end{center}\vspace{2mm}\end{minipage}\end{Sbox}\fbox{\TheSbox}\vspace{3mm}

We now define a function $v$ that satisfies \Hpp\ and a measure $\mu_0\in\Prr$, and we show two solutions of \r{e-cauchy}. The example is based on the classical example of non-uniqueness of the solution for the Cauchy problem $\dot x=\sqrt{x},\ x(0)=0$.

\begin{figure}[htb]
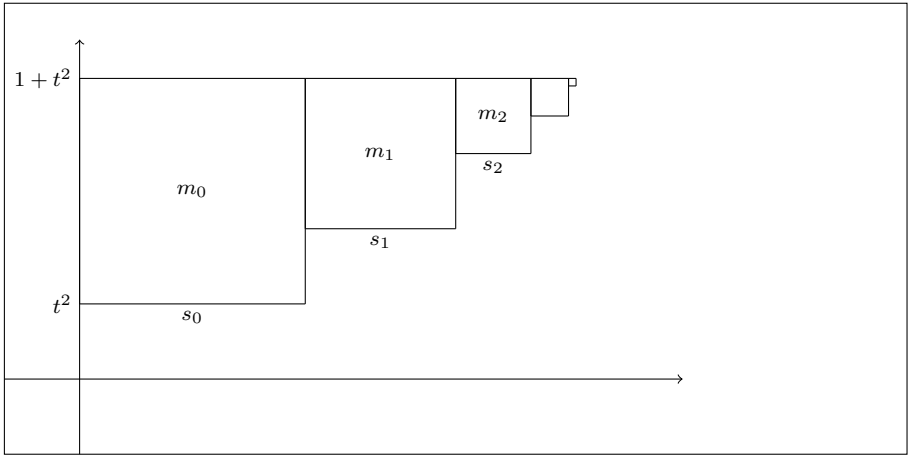

\begin{center} \begin{pgfpictureboxed}{0cm}{0cm}{12 cm}{6 cm}

\begin{pgfscope}
\pgfsetendarrow{\pgfarrowto}
\pgfline{\pgfxy(1,0)}{\pgfxy(1,5.5)}
\pgfline{\pgfxy(0,1)}{\pgfxy(9,1)}
\end{pgfscope}
\quadrato{1}{2}{3}
\quadrato{4}{3}{2}
\quadrato{6}{4}{1}
\quadrato{7}{4.5}{.5}
\quadrato{7.5}{4.9}{.1}

\pgfputat{\pgfxy(2.5,1.9)}{\pgfbox[center,top]{$s_0$}}
\pgfputat{\pgfxy(5,2.9)}{\pgfbox[center,top]{$s_1$}}
\pgfputat{\pgfxy(6.5,3.9)}{\pgfbox[center,top]{$s_2$}}

\pgfputat{\pgfxy(2.5,3.5)}{\pgfbox[center,center]{$m_0$}}
\pgfputat{\pgfxy(5,4)}{\pgfbox[center,center]{$m_1$}}
\pgfputat{\pgfxy(6.5,4.5)}{\pgfbox[center,center]{$m_2$}}

\pgfputat{\pgfxy(.9,2)}{\pgfbox[right,center]{${t}^{2}$}}
\pgfputat{\pgfxy(.9,5)}{\pgfbox[right,center]{$1+{t}^{2}$}}
\end{pgfpictureboxed}

\caption{\Hpp\ does not guarantee uniqueness of the solution.}
\label{fig:controL1}
\end{center}
\end{figure}

We first fix the dimension of the space $n=2$ and the final time $T=1$. We define a family of measures $\nu_t\in\Prr$ as shown in Figure \ref{fig:controL1} (dimensions are not respected).\\

More precisely, $\nu_t$ is defined as follows. First define a sequence of squares $Q^i_t$, all with sides parallel to axes $x$ and $y$. Denote with $s_i$ the length of the side of the square $Q^i_t$. All the squares have the upper side on the line $y=1+t^{2}$. The square $Q^0_t$ has the left side on the $y$-axis. For $i>0$, the left side of the square $Q^i_t$ is contained in the right side of the square $Q^{i-1}_t$. Define $$\nu_t:=\sum_{i=0}^\infty  m_i \chi_{Q^i_t} \lam,$$ where $m_i$ are positive real numbers to be chosen, and $\lam$ is the Lebesgue measure.

We now choose $s_i:=4^{-i}$ and $m_i=\frac{1}{2} 8^i$. It is evident that $\nu_t$ is positive and absolutely continuous with respect to the Lebesgue measure. It is easy to prove that $\nu_t$ has bounded support in a rectangle of sides $s_0$ and $\sum_{i=0}^\infty s_i=\frac43$. We also that $\nu_t(\R^n)=\sum_{i=0}^\infty m_i s_i^2=1$.

We now define $v\Pq{\nu_t}:=(0,2t)$. We first prove that $v$ satisfies \Hpp\ for $\nu_t$ defined before. We then extend $v$ to the whole space $\Pr$ so that it satisfies \Hpp\ everywhere. For each $\nu_t$, $v\Pq{\nu_t}$ is a constant vector field, thus $v$ is uniformly Lipschitz with $L=0$ and uniformly bounded with $M=2T=2$. We now prove that $v$ is Lipschitz with respect to $L^1$ distance. Let $t,s\in\Pq{0,1}$ and assume $t>s$ with no loss of generality. Then, for all $x\in\R^n$ we have
$|v\Pq{\nu_t}(x)-v\Pq{\nu_s}(x)|=2(t-s).$
Take now $n$ such that $s_n<t^{2}-s^{2}\leq s_{n-1}$. Remark that such a $s_n$ always exists, since $s_n\searrow 0$. Then, for all $i\geq n$ the square $Q^i_t$ is disjoint with respect to $Q^i_s$. Moreover, all $Q^i_t$ are always disjoint with respect to $Q^j_s$ for $j\neq i$. As a consequence,
\bqn
\|\nu_t-\nu_s\|_{L^1}\geq \sum_{i=n}^\infty 2 \nu_t(Q^i_t)=\sum_{i=n}^\infty 2 m_i s_i^2=2 \cdot 2^{-n},
\eqnn
hence
\bqn
\|v\Pq{\nu_t}-v\Pq{\nu_s}\|_{C^0}&=&2\Pt{t-s}\leq 2\sqrt{t^2-s^2}\leq 2 \sqrt{s_{n-1}} =2\cdot 2^{-n+1} \leq 2 \|\nu_t-\nu_s\|_{L^1}.
\eqnn
Thus $v$ is Lipschitz with respect to $L^1$ distance, with $K=2$.

We have just proved that $v$ satisfies \Hpp\ for the $\nu_t$. We now observe that we have the two following solutions of \r{e-cauchy} with $\mu_0=\nu_0$. The first is the constant solution $\mu_t=\nu_0$ for all $t\in\Pq{0,1}$, the second is $\mu_t=\nu_t$. It is easy to prove that the first is a solution, since $v\Pq{\nu_0}=(0,0)$, thus constant solutions are solutions. For the second, we have to prove that $\nu_t$ satisfies \r{e-cauchy} with $v\Pq{\nu_t}=(0,2t)$. It is equivalent to define a non-autonomous vector field $w(t,x):=(0,2t)$, consider the flow that it generates $\Phi^w_t(x,y)=(x,y+t^2)$ and observe that $\Phi^w_t\#\nu_0=\nu_t$.

It is easy to observe that $v$ does not satisfy \Hp, for all $p\geq 1$. Indeed, $W_p(\nu_t,\nu_s)=t^2-s^2$, thus it is impossible to find a finite $K'$ such that, for all $t,s$, it holds $$\|v\Pq{\nu_t}-v\Pq{\nu_s}\|_{C^0}=2(t-s)\leq K'(t^2-s^2)=K' W_p(\nu_t,\nu_s).$$

We now extend the definition of $v$ for any $\mu\in\Pr$. Consider the rectangle $R:=[0,4/3]\times[1,2]$ and define the following functional
\mmfunz{F}{\Pr}{[0,+\infty)}{\mu}{F(\mu):=\mu\Pt{R}.}
Observe that $F$ is well defined and finite, since $F(\mu)\leq \mu(\R^2)=1$. We now define the function $\fz{f}{\Pq{0,1}}{\R}$ as follows
$$f(F(\nu_t)):=2t,$$
where the $\nu_t$ are defined above. It is easy to prove that this function is well defined, and defined on the whole domain. We now prove that $f$ is a Lipschitz function. We take $0<\xi_1<\xi_2\leq 1$ and prove that $|f(\xi_2)-f(\xi_1)|\leq K |\xi_2-\xi_1|$. We have
\bqn
\frac{|f(\xi_2)-f(\xi_1)|}{|\xi_2-\xi_1|}=\frac{2|t_2-t_1|}{|F(\nu_{t_2})-F(\nu_{t_1})|},
\eqnn
where $t_i$ is the unique value such that $F(\nu_{t_i})=\xi_i$. We first study the particular case $s_{n+1}\leq t_1^2 <t_2^2 \leq s_{n}$, where $s_i=4^{-i}$ has been defined above. We have that all squares $Q^i_{t_k}$ with $i\geq n+1$, $k=1,2$, are completely contained in $R$, since $1+t_k^2-s_m\geq 1$. On the other side, each square $Q^i_{t_k}$ with $i\leq n$ has an intersection with $R$ that is given by the rectangle of horizontal side $s_i$ and vertical side $t_k^2$. As a consequence, one has 
\bqn
|F(\nu_{t_2})&-&F(\nu_{t_1})|=\Pabs{\sum_{i=0}^{n} \nu_{t_2}(Q^i_{t_2}\cap R) - \sum_{i=0}^{n} \nu_{t_1}(Q^i_{t_1}\cap R)}=\nn
&&=\Pabs{\sum_{i=0}^{n} m_i s_i t_2^2 - \sum_{i=0}^{n} m_i s_i t_1^2}=|t_2-t_1|\,|t_2+t_1| \sum_{i=0}^{n} m_i s_i\geq 2 |t_2-t_1|.
\eqnn
 We thus have
\bqn
\frac{|f(\xi_2)-f(\xi_1)|}{|\xi_2-\xi_1|}\leq \frac{2|t_2-t_1|}{|t_2-t_1| 2^{-1}}=4,
\eqnn
for $\xi_1,\xi_2$ satisfying $s_{n+1}\leq t_1^2<t_2^2\leq s_{n}$ for some $n\in \N$. Since the Lipschitz constant $4$ does not depend on $n$, we can easily pass to the general case via the triangular inequality. One can observe that the Lipschitzianity is also verified in $0$, since $f$ is continuous and the Lipschitz constant $4$ for $\xi>0$ does not depend on $\xi$ itself.

We now define $v\Pq{\mu}(x):=(0,f(F(\mu)))$. By construction, we have $v\Pq{\nu_t}(x)=(0,2t)$, i.e. this definition of $v$ coincides with the definition given above, that gives non-uniqueness. We now prove that such definition of $v$ satisfies \Hpp. The only non-straightforward point is to prove that $v$ is Lipschitz with respect to the $L^1$ distance. Take two distinct measures $\mu,\nu$. If they satisfy $F(\mu)=F(\nu)$, then 
$\|v\Pq{\mu}-v\Pq{\nu}\|_{C^0}=0<\|\mu-\nu\|_{L^1}$. Otherwise, observe that $
\|\mu-\nu\|_{L^1}\geq \int_{R} |\mu(x)-\nu(x)|\,d\lam(x)\geq|F(\mu)-F(\nu)|.$  Thus $\frac{\|v\Pq{\mu}-v\Pq{\nu}\|_{C^0}}{\|\mu-\nu\|_{L^1}}\leq \frac{|f(F(\mu))-f(F(\nu))|}{|F(\mu)-F(\nu)|}$. Hence, $v$ is uniformly Lipschitz, since $f$ is Lipschitz.

\brem One can ask if condition \Hpp\ is sufficient to guarantee at least existence of a solution for \r{e-cauchy}. This is, to our knowledge, an open question. Nevertheless, the previous discussion is sufficient to show that $L^1$ distance is not a good choice for the study of \r{e-cauchy} in the case of $v$ depending on the measure itself.\\
\erem

%We now prove that hypotheses \Hp\ with $p=1$ imply hypotheses \Hpp. Indeed, consider two distribution $\mu,\nu$ and remark that $W_1(\mu,\nu)=W_1(\mu-\Pq{\mu-\nu}_+,\nu-\Pq{\nu-\mu}_+)$, that means that the $W_1$ distance can be computed by removing the shared mass. This is a particular case of invariance of $W_1$ distance under mass subtraction, see \cite[Cor. 1.16]{villani}. Now assume that $\mu,\nu$ have no shared mass. Thus $\|\mu-\nu\|_{L_1}=2\mu(\R^2)$. We also have
%$W_1(\mu,\nu)=\int |x-T(x)|\, d\mu(x)$. Since $x\in\supp(\mu)$ and $T(x)\in \supp(\nu)$, we have
%\bqn
%W_1(\mu,\nu)\leq \diam(\supp(\mu)\cup\supp(\nu)) \int\,d\mu(x)=\frac{\diam(\supp(\mu)\cup\supp(\nu))}{2} \|\mu-\nu\|_{L^1}.
%\eqnn
%The estimate is finite, since we deal with $\mu,\nu\in\Mu$, and in particular they have compact support. We don't have a global constant $K$ such that $W_1(\mu,\nu)\leq K \|\mu-\nu\|_{L^1}$, but for all bounded set $B\subset \Mu$ with respect to $

%, while the converse is not true. Indeed, take a solution $\mu_t$ of \r{e-cauchy}. For the moment, we don't know if it is unique. let $D_t$ be the diameter of the support of $\mu_t$. We know that $D_0 < \infty$ and $D_t\leq D_0+2Mt$. In particular, $D_t\leq D_0+2MT$ for all $t\in\Pq{0,T}$.

\noindent \b{Acknowledgments}: This work was conducted during a  visit of F. Rossi to Rutgers University, Camden, NJ, USA. He thanks the institution for its hospitality.

\noindent The authors thank the anonymous reviewer for remarks and the suggested bibliography.

\end{document}